\documentclass[10pt]{amsart}
\usepackage[dvips]{epsfig}
\usepackage{amscd}
\usepackage{amssymb}
\usepackage{amsthm}
\usepackage{amsmath}
\usepackage{latexsym}

\usepackage{upref}
\usepackage{hyperref}

\theoremstyle{plain}
\newtheorem{thm}{Theorem}[section]
\theoremstyle{plain}
\newtheorem{lem}[thm]{Lemma}
\newtheorem{prop}[thm]{Proposition}

\theoremstyle{definition}
\newtheorem{defi}{Definition}[section]
\newtheorem{rem}{Remark}

{%
\setcounter{enumi}{0}

\begin{enumerate}}%
{\end{enumerate} }
{%
\setcounter{enumi}{0}

\begin{enumerate}}%
{\end{enumerate} }

\newcommand{\norm}[1]{\ensuremath{\left\|#1\right\|}}
\newcommand{\abs}[1]{\ensuremath{\left|#1\right|}}
\newcommand{\Om}{\ensuremath{\Omega}}
\newcommand{\disp}{\ensuremath{\displaystyle}}

\newcommand{\TOm}{\ensuremath{Q_T}}
\newcommand{\pTOm}{\ensuremath{\partial \Om}\times (0,T)}
\newcommand{\pOm}{\ensuremath{\partial \Om}}

\newcommand{\bu}{\ensuremath{{\bf u}}}
\newcommand{\bv}{\ensuremath{{\bf v}}}

\newcommand{\D}{\ensuremath{\mathcal{D}}}

\newcommand{\R}{\ensuremath{\Bbb{R}}}
\newcommand{\N}{\ensuremath{\Bbb{N}}}

\newcommand{\Grad}{\mathrm{\nabla}}

\newcommand{\sn}{\sum_{n=0}^{N-1}}
\newcommand{\si}{\sum_{i=1}^3}

\newcommand{\pt}{\ensuremath{\partial_t}}

\newcommand{\dx}{\ensuremath{\, dx}}

\newcommand{\dt}{\ensuremath{\, dt}}
\newcommand{\ds}{\ensuremath{\, ds}}

\numberwithin{equation}{section} \allowdisplaybreaks

\title[Finite volume schemes for nonlocal reaction-diffusion systems]
{Convergence of a finite volume scheme for nonlocal reaction-diffusion
systems modelling an epidemic disease}

\date{\today}
\author[Mostafa Bendahmane]{Mostafa Bendahmane}
\address[Mostafa Bendahmane]{\newline
         CI$^2$MA and Departamento de Ingenier\'{\i}a Matem\'atica\newline
         Universidad de Concepci\'on\newline
         Casilla 160-C, Concepci\'on, Chile}
\email[]{mostafab@ing-mat.udec.cl}

\author[M. A. Sep\'ulveda C.]{Mauricio A. Sep\'ulveda}
\address[Mauricio A. Sepulveda]{\newline
         CI$^2$MA and Departamento de Ingenier\'{\i}a Matem\'atica\newline
         Universidad de Concepci\'on\newline
         Casilla 160-C, Concepci\'on, Chile}
\email[]{mauricio@ing-mat.udec.cl}


\thanks{Partially supported by CMM, Universidad de Chile, and CI$^2$MA, Universidad de Concepci\'on.}

\begin{document}

\begin{abstract}
We analyze a finite volume scheme for nonlocal SIR model, which 
is a nonlocal reaction-diffusion system modeling an epidemic disease. 
We establish existence solutions to the finite volume scheme, and show that it converges to a weak solution.
The convergence proof is based on deriving 
series of a priori estimates and using a general $L^p$ compactness criterion.
\end{abstract}

\maketitle


\section{Introduction}

We consider a mathematical model describing an epidemic disease in a physical domain 
$\Om\subset \mathbb{R}^d$ ($d=1,2,3$) over a time span $(0,T)$, $T>0$.
In this model we consider the propagation of an epidemic disease in a 
simple population $p=u_1+u_2+u_3$, where $u_1=u_1(t,x)$,  $u_2=u_2(t,x)$ and
$u_3=u_3(t,x)$ are the respective densities of susceptible (those who can
catch the disease), infected (those who have the disease and
can transmit it) and of recover individuals (those who have been
exposed to the disease and will become infective after the lapse of an incubation period) at time $t$ and location $x$.
A prototype of a nonlinear system that governs the spreading of a nonlocal SIR of epidemics with
through a population in a spatial domain is the following nonlocal reaction-diffusion system:
\begin{equation}
	\label{S1}
	\begin{cases}
		&\displaystyle \pt u_1 
		-a_1\Bigl(\int_\Om u_1\dx\Bigl) \Delta u_1= -\sigma(u_1,u_2,u_3)-\mu u_1,\\ 
		&\displaystyle \pt u_2 
		-a_2\Bigl(\int_\Om u_2\dx\Bigl) \Delta u_2 = \sigma(u_1,u_2,u_3)-\gamma u_2-\mu u_2,\\ 
		&\displaystyle \pt u_3 
		-a_3\Bigl(\int_\Om u_3\dx\Bigl)\Delta u_3= \gamma u_2,
	\end{cases}	
\end{equation}
in $\TOm$, where $\TOm$ denotes the time-space cylinder $(0,T)\times \Om$. 
We complete the system \eqref{S1} with Neumann boundary conditions:
\begin{equation}
	\label{S2}
	a_i\Bigl(\int_\Om u_i\dx\Bigl)\Grad u_i\cdot \eta=0 \quad \text{on $\pTOm$}, \qquad  i=1,2,3,
\end{equation} 
where $\eta$ denotes the outer unit normal to the boundary $\pOm$ of $\Omega$, and with initial data:
\begin{equation}
	\label{S3}
	u_i(0,x)=u_{i,0}(x), \quad x\in \Om, \qquad  i=1,2,3.
\end{equation}
One of the simplest SIR models is the Kermack-McKendrick model \cite{Kermack}
which consists in a systems of $3\times 3$ EDOs which is similar to the
model above, but without diffusion terms.
The classical epidemic SIR model describes the infection and recovery process 
in terms of three ordinary differential equations has been studied by many 
researchers \cite{diekmann,wang} and the reference cited therein. 

Hence the description of the classical SIR model
in the system \eqref{S1}, $\frac{1}{\gamma}$ is
the length of latency period or duration of the exposed stage, and
$\mu$ the natural mortality rate. 
The incidence term $\sigma$ take the following form:
\begin{equation}\label{def-sigma-1}
\sigma(u,v,w)=\alpha \frac{uv}{u+v+w} \text{ for some } \alpha>0,
\end{equation}
which coincides with the classical model  $\sigma(u,v,w)=\alpha {uv}$  when the total population $P(t)$ remains constant.
In fact, the well known SIR model which appears generally in literature is renormalized
which we do not suppose here.
For technical reasons, we need to extend the function
$\sigma(u,v,w)$ so that it becomes defined for all $(u,v,w)\in
{\R}\times{\R}\times{\R}$. We do this by setting
$$
\sigma(u,v,w)=
\left\{
\begin{array}{lll}
\sigma(u^+,v^+,w^+)&\mbox{ if } &(u,v,w)\not= (0,0,0),\\
0&\mbox{ otherwise }.
\end{array}
\right.
$$

In this work, the diffusion  rates $a_1 > 0$, $a_2 > 0$ and $a_3 > 0$
are supposed to depend to the whole of each populations in the domain rather than on the
local density, i. e. moves are guided by considering the global state
of the medium.
We assume that
$a_i:\R\to\R$ is a continuous function satisfying: there exist constants $M_i,C>0$ such that
\begin{equation}
	\label{ass-ai}
M_i \leq a_i\text{ and } \abs{a_i(I_1)-a_i(I_2)}\leq C
\abs{I_1-I_2}\text{ for all $I_1,I_2 \in \R$, for $i=1,2,3$}.
\end{equation}
Such  equations with nolocal diffusion terms 
has already been studied from a theoretical point of view by several authors.
First, in 1997, M. Chipot and B. Lovat
\cite{chipot} studied the existence and uniqueness of the solutions
for a scalar parabolic equation with a nonlocal diffusion term.
Existence-uniqueness and long time behaviour for other class of nonlocal nonlinear
parabolic equations and systems are studied in \cite{ackleh,raposo}.
Liu and Jin made some experimental simulations  in \cite{liu} in order to observe
spacial patterns in an epidemic model with constant removal rate of the infective. 
\\

Before we define our finite volume scheme, let us state 
a relevant definition of a weak solution for the nonlocal SIR model.

\begin{defi} \label{def1}\label{def1}
A weak solution of \eqref{S1}-\eqref{S3} is a triple 
$\bu=(u_1,u_2,u_3)$ of functions such that $u_1,u_2,u_3 \in L^{2}(0,T;H^1(\Om))$, 
\begin{equation}\label{eq1:def1}
     \begin{split}
     &\displaystyle-\int_{\Om} u_{1,0}(x) \varphi_1(0,x)\dx
     -\iint_{Q_T} u_1 \, \pt \varphi_1\dx\dt+\iint_{Q_T} a_1\Bigl(\int_\Om u_1\dx\Bigl)
     \Grad u_1\cdot \Grad \varphi_1~\dx\dt\\
     &\qquad \qquad =-\iint_{Q_T} (\sigma(u_1,u_2,u_3)+\mu u_1)\varphi_2~\dx\dt,
    \end{split}
\end{equation}
\begin{equation}\label{eq2:def1}
     \begin{split}
     &\displaystyle-\int_{\Om} u_{2,0}(x) \varphi_2(0,x)\dx
     -\iint_{Q_T} u_2 \, \pt \varphi_2\dx\dt+\iint_{Q_T} a_2\Bigl(\int_\Om u_2\dx\Bigl)
     \Grad u_2\cdot \Grad \varphi_2~\dx\dt\\
     &\qquad \qquad =\iint_{Q_T} (\sigma(u_1,u_2,u_3)-(\gamma+\mu) u_2)\varphi_2~\dx\dt,
    \end{split}
\end{equation}
\begin{equation}\label{eq3:def1}
     \begin{split}
     &\displaystyle-\int_{\Om} u_{3,0}(x) \varphi_3(0,x)\dx
     -\iint_{Q_T} u_3 \, \pt \varphi_3\dx\dt+\iint_{Q_T} a_3\Bigl(\int_\Om u_3\dx\Bigl)
     \Grad u_3\cdot \Grad \varphi_3~\dx\dt\\
     &\qquad \qquad =\iint_{Q_T} \gamma u_2 \varphi_3~\dx\dt,
    \end{split}
\end{equation}
for all $\varphi_1,\varphi_2,\varphi_3 \in \D([0,T)\times \overline{\Om})$.
\end{defi}

\begin{rem}
Note that we can easily check that Definition \ref{def1} makes sense. 
Furthermore, observe that Definition \ref{def1} implies that $\pt u_i$ 
belongs to $L^2\left(0,T;(H^1(\Om))'\right)$, so 
that $u_i\in C([0,T];L^2(\Om))$ for $i=1,2,3$.
\end{rem}

\begin{rem}
A classical way to prove the existence of weak solutions in the sense  of \eqref{eq1:def1}-\eqref{eq3:def1},
is to use Faedo-Galerkin method like the system studied in  \cite{ackleh} or in \cite{raposo}.
On the other hand, the proof here, of convergence of the numerical scheme,
implies the existence of weak solutions of \eqref{S1}-\eqref{S3}. Additionally, a proof of uniqueness 
of the weak solution is given 
in the appendix.
\end{rem}

Following \cite{Ey-Gal-Her:book}, we now give a precise definition of the finite volume scheme 
for the nonlocal SIR model. Let $\Omega$ be an open bounded polygonal connected subset of
$\mathbb{R}^3$ with boundary $\partial \Omega$. Let $\Omega_R$ be an admissible mesh of 
the domain $\Om$ consisting of open and convex polygons called control 
volumes with maximum size (diameter) $h$.
For all $K \in \Omega_R$, let by $x_K$ denote the center of $K$, $N(K)$ the set of the
neighbors of $K$ i.e. the set of cells of $\Omega_R$ which have a common
interface with $K$, by $N_{\text{int}}(K)$ the set of the neighbors of $K$
located in the interior of $\Omega_R$, by $N_{\text{ext}}(K)$ the set of
edges of $K$ on the boundary $\partial \Om$. Furthermore, for all $L \in N(K)$
denote by $d(K,L)$ the distance between $x_K$ and $x_{L}$, by $\sigma_{K,L}$
the interface between $K$ and $L$, by $\eta_{K,L}$ the unit normal vector
to $\sigma_{K,L}$ outward to $K$. For all $K \in \Omega_R$, we denote by
$m(K)$ the measure of $K$. The admissibility of $\Omega_R$ implies that
$\overline{\Om}=\cup_{K\in \Omega_R} \overline{K}$, $K\cap L=\emptyset$
if $K,L\in \Omega_R$ and $K \ne L$, and there exist a finite sequence
of points $(x_{K})_{K\in \Omega_R}$ and the straight line $\overline{x_{K}x_{L}}$
is orthogonal to the edge $\sigma_{K,L}$. 
We also need some regularity on the mesh:
$$
\min_{K\in \Omega_R,L \in N(K)}
\frac{d(K,L)}{\text{diam}(K)}\ge \alpha
$$
for some $\alpha>0$. 


We denote by $H_h(\Om) \subset L^2(\Om)$ the space of functions
which are piecewise constant on each control volume $K \in \Omega_R$. For all
$u_h \in H_h(\Om)$ and for all $K \in \Omega_R$, we denote by $u_K$ the constant value
of $u_h$ in $K$. For $(u_h,v_h)\in (H_h(\Om))^2$, we define the following inner product:
$$
\left \langle u_h,v_h \right \rangle_{H_h}= \frac{1}{2}\sum_{K \in \Omega_R }\sum_{L \in N(K) }
\frac{m(\sigma_{K,L})}{d(K,L)}(u_{L}-u_{K})(v_{L}-v_{K}),
$$
corresponding to Neumann boundary conditions. We define a norm in $H_h(\Om)$ by
$$
\norm{u_h}_{H_h(\Om)}=(\left \langle u_h,u_h \right \rangle_{H_h})^{1/2}.
$$
Finally, we define $L_h(\Om) \subset L^2(\Om)$ the space of functions
which are piecewise constant on each control volume $K \in \Omega_R$ with the associated
norm
$$
\left (u_h,v_h \right )_{L_h(\Om)}= \sum_{K \in \Omega_R }
m({K})u_{K} v_K,\qquad \norm{u_h}^2_{L_h(\Om)}
=\sum_{K \in \Omega_R }m({K})\abs{u_{K}}^2,
$$
for $(u_h,v_h)\in (L_h(\Om))^2$. 

Next, we let $K \in \Omega_R$ and $L \in N(K)$ with common vertexes
$(a_{\ell,K,L})_{1\le \ell\le I}$ with $I \in \N^{\star}$. Next let
$T_{K,L}$ (respectively $T^{\text{ext}}_{K,\sigma}$ for $\sigma\in N_{\text{ext}}(K)$) be the
open and convex polygon with vertexes $(x_K,x_L)$ ($x_K$ respectively) and
$(a_{\ell,K,L})_{1\le \ell\le I}$. Observe that
$$
\displaystyle \Om=\cup_{K\in \Omega_R}
\Biggl(\Bigl(\cup_{L\in N(K)}\overline{T}_{K,L}\Bigl)
\cup \Bigl(\cup_{\sigma \in N_{\text{ext}}(K)}\overline{T}^{\text{ext}}_{K,\sigma}\Bigl)\Biggl)
$$

The approximation $\Grad_h u_{h}$ of $\Grad u$ is defined by
$$
\Grad_h u_{h}(x)=\begin{cases}
\frac{m(\sigma_{K,L})}{d(K,L)}(u_{L}-u_{K})\eta_{K,L}&\text{ if $x \in T_{K,L}$},\\
0&\text{ if $x \in T^{\text{ext}}_{K,\sigma}$},
\end{cases}
$$
for all $K \in \Omega_R$.

The next goal is to discretize the problem \eqref{S1}-\eqref{S3}.
We denote by $\mathcal{D}$ an admissible discretization of $Q_T$, which consists of
an admissible mesh of $\Om$, a time step $\Delta t>0$, and a positive number $N$ chosen 
as the smallest integer such that $N\Delta t\ge T$.  
We set $t^n=n \Delta t$ for $n\in  \{0,\ldots,N\}$. 

We approximate the nonlocal SIR model in the following 
way: Determine vectors $(u_{i,K}^{n})_{K \in \Omega_R}$ for $i=1,2,3$,
such that for all $K \in \Omega_R$ and $n \in  \{0,\ldots,N-1\}$ 
\begin{equation}\label{prob:init}
u_K^{i,0}=\frac{1}{m(K)} \int_{K} u_{i,0}(x) \dx,\quad i=1,2,3,
\end{equation}
\begin{equation}\label{S1-discr}
     \begin{split}
     &m(K)\frac{u^{n+1}_{1,K}-u^{n}_{1,K}}{\Delta t}
	 -a_1\Bigl(\sum_{K_0 \in \Om_h}u^{n}_{1,K_0}\Bigl)\sum_{L \in N(K) }
     \frac{m(\sigma_{K,L})}{d(K,L)}(u^{n+1}_{1,L}-u^{n+1}_{1,K})
     \\ &\qquad \qquad \qquad
     +m(K)\left(\sigma(u^{n+1,+}_{1,K},u^{n+1,+}_{2,K},u^{n+1,+}_{3,K})+\mu u^{n+1}_{1,K}\right)=0,
     \end{split}
\end{equation}
\begin{equation}\label{S2-discr}
     \begin{split}
     &m(K)\frac{u^{n+1}_{2,K}-u^{n}_{2,K}}{\Delta t}
	 -a_2\Bigl(\sum_{K_0 \in \Om_h}u^{n}_{2,K_0}\Bigl)\sum_{L \in N(K) }
     \frac{m(\sigma_{K,L})}{d(K,L)}(u^{n+1}_{2,L}-u^{n+1}_{2,K})
     \\ &\qquad \qquad 
     -m(K)\left(\sigma(u^{n,+}_{1,K},u^{n+1,+}_{2,K},u^{n+1,+}_{3,K})-(\gamma+\mu) u^{n+1}_{2,K}\right)=0,
     \end{split}
\end{equation}
\begin{equation}\label{S3-discr}
     \begin{split}
     &m(K)\frac{u^{n+1}_{3,K}-u^{n}_{3,K}}{\Delta t}
	 -a_3\Bigl(\sum_{K_0 \in \Om_h}u^{n}_{3,K_0}\Bigl)\sum_{L \in N(K) }
     \frac{m(\sigma_{K,L})}{d(K,L)}(u^{n+1}_{3,L}-u^{n+1}_{3,K})
     \\ &\qquad \qquad \qquad\qquad \qquad\qquad \qquad\qquad \qquad\qquad
     -m(K)\gamma u^{n}_{2,K}=0.
     \end{split}
\end{equation}



To simplify the notation, it will always be understood that when $h$ is sent to zero then so 
is $\Delta t$, thereby assuming (without loss of generality) a functional relationship 
between the spatial and temporal discretization parameters. This is not a real restriction on 
the time step, but it allows us to write ``$u_{i,h}$" instead of
``$u_{i,h,\Delta t}$" for $i=1,2,3$ , ``$h\to 0$" 
instead of ``$h,\Delta t\to 0$", and so forth.

For the sake of analysis, we introduce the ``piecewise constant" functions 
\begin{equation}\label{prob:general}
	u_{i,h}(t,x)=u_{i,K}^{n+1},~~i=1,2,3,
\end{equation}
for all $(t,x) \in (n \Delta t,(n+1)\Delta t]\times K$, with $K \in \Omega_R$ 
and $n \in  \{0,\ldots,N-1\}$. 
To simplify the notation, let us write 
$\bu_h$ for the vector $(u_{1,h},u_{2,h},u_{3,h})$.
Our main result is
\begin{thm} \label{theo1}
Assume $u_{i,0}\in L^2(\Om)$ for $i=1,2,3$. 
Then the finite volume solution $\bu_h$, generated 
by \eqref{prob:init} and \eqref{S1-discr}-\eqref{S3-discr}, converges 
along a subsequence to $\bu=(u_1,u_2,u_3)$ as $h\to 0$, 
where $\bu$ is a weak solution of \eqref{S1}-\eqref{S3}. 
The convergence is understood in the following sense:  
\begin{equation*}
\begin{split}
&\text{$u_{i,h}\to u_i$ strongly in $L^2(Q_T)$ and a.e. in $Q_T$},\\
&\text{$\Grad_h u_{i,h}\to \Grad u_{i}$ weakly in $(L^2(Q_T))^3$},\\
&\text{$\sigma(u_{1,h},u_{2,h},u_{3,h})\to \sigma(u_{1},u_{2},u_{3})$
strongly in $L^1(Q_T)$},
\end{split}
\end{equation*}
for $i=1,2,3$.
\end{thm}

The remaining part of this paper is organized as follows.
The proof of Theorem \ref{theo1} is divided into Section \ref{sec-exist-est}
(existence of the scheme),
Section \ref{sec:basic-apriori} (basic a priori estimates), Section \ref{sec-compactness}
(space and time translation estimates), and Section \ref{sec:conv} (convergence to a weak solution).
In section \ref{num-ex} we give some numerical examples. Finally in Appendix we prove
the uniqueness of the solution using duality techniques.

\section{Existence of the finite volume scheme}\label{sec-exist-est}

The existence of a solution to the finite volume scheme 
\eqref{prob:init}-\eqref{S3-discr} will be obtained with the help of the following lemma
proved in \cite{Lions:1969} and \cite{Tem:2001}.

\begin{lem}\label{lem:exist-classic}
Let $\mathcal{A}$ be a finite dimensional Hilbert space with scalar product $[\cdot,\cdot]$
and norm $\norm{\cdot}$, and let $\mathcal{P}$ be a continuous mapping from $\mathcal{A}$
into itself such that
$$
[\mathcal{P}(\xi),\xi]>0 \text{ for } \norm{\xi} =r>0.
$$
Then there exists $\xi \in \mathcal{A}$ with $\norm{\xi} \le r$ such that
$$
\mathcal{P}(\xi)=0.
$$
\end{lem}

The existence for the finite volume scheme is given in

\begin{prop}
\label{prop:exist-fv}
Let $\mathcal{D}$ be an admissible discretization of $Q_{T}$. 
Then the problem \eqref{prob:init}-\eqref{S3-discr} admits at least one solution
$(u^{n}_{1,K},u^{n}_{2,K},u^{n}_{3,K})_{(K,n) \in \Omega_R\times  \{0,\ldots,N\}}$.
\end{prop}

\begin{proof}

First we introduce the Hilbert spaces
$$E_h=(H_h(\Om)\cap L_h(\Om))^3, 
$$ 
under the norm
$$
\norm{\bu_h}^2_{E_h}:=\si \norm{u_{i,h}}^2_{H_h(\Om)}+
\si \sum_{K \in \Omega_R }m(K)\abs{u_{i,K}}^2,
$$
where $\bu_h=(u_{1,h},u_{2,h},u_{3,h})$.
Let $\Phi_h=(\varphi_{1,h},\varphi_{2,h},\varphi_{3,h})\in E_h$ and define the discrete bilinear forms
$$
T_{h}(\bu_{h}^{n},\Phi_{h})=
\si \Bigl(u^{n}_{i,h},\varphi_{i,h}\Bigl),
$$
\begin{equation*}\begin{split}
b_{h}(\bu_{h}^{n+1},\Phi_{h})=&\sum_{K \in \Omega_R }m(K)
\Bigl(\sigma(u^{n+1,+}_{1,K},u^{n+1,+}_{2,K},u^{n+1,+}_{3,K})+\mu u^{n+1}_{1,K}\Bigl)\varphi_{1,K}\\
&-\sum_{K \in \Omega_R }m(K) 
\Bigl(\sigma(u^{n,+}_{1,K},u^{n+1,+}_{2,K},u^{n+1,+}_{3,K})-(\gamma+\mu) u^{n+1}_{2,K}\Bigl)
\varphi_{2,K},\\
&-\gamma \sum_{K \in \Omega_R }m(K) u^{n}_{2,K}\varphi_{3,K},
\end{split}\end{equation*}
and
$$
a_h(\bu_{h}^{n+1},\Phi_{h})=\si
\frac{\displaystyle a_i\Bigl(\sum_{K_0 \in \Om_h}u^{n}_{i,K_0}\Bigl)}{2}
\sum_{K \in \Omega_R }\sum_{L \in N(K) }
\frac{m(\sigma_{K,L})}{d(K,L)}(u^{n+1}_{i,L}-u^{n+1}_{i,K})
(\varphi^{n+1}_{i,L}-\varphi^{n+1}_{i,K}).
$$
Multiplying \eqref{S1-discr}, \eqref{S2-discr} and \eqref{S3-discr}
by $\varphi_{1,K}$, $\varphi_{2,K}$, $\varphi_{3,K}$, respectively, we get 
the equation
\begin{equation*}
\frac{1}{\Delta t}\Bigl(T_{h}(\bu_{h}^{n+1},\Phi_{h})-T_{h}(\bu_{h}^{n},\Phi_{h})\Bigl)\\
+a_h(\bu_{h}^{n+1},\Phi_{h})+b_{h}(\bu_{h}^{n+1},\Phi_{h})=0.
\end{equation*}
Now we apply the Lemma \ref{lem:exist-classic} for proving the existence of
$\bu_{h}^{n+1}$ for all $K \in \Omega_R$ and $n \in  \{0,\ldots,N\}$.
We define the mapping $\mathcal{P}$ from $E_h$ into itself 
\begin{equation*}
\begin{split}
[\mathcal{P}(\bu_{h}^{n+1}),\Phi_{h}]=&
\frac{1}{\Delta t}(T_{h}(\bu_{h}^{n+1},\Phi_{h})-T_{h}(\bu_{h}^{n},\Phi_{h}))\\
&\qquad \qquad \qquad +a_h(\bu_{h}^{n+1},\Phi_{h})+b_{h}(\bu_{h}^{n+1},\Phi_{h}),
\end{split}
\end{equation*}
for all $\Phi_{h} \in E_h$.
Note that it is easy to obtain from the discrete H\"{o}lder inequality the following 
bounds:
\begin{align*}
a_h(\bu_{h},\bv_h )& \le C \norm{\bu_{h}}_{E_h}  \norm{\bv_h}_{E_h},
\\ T_h(\bu_{h},\bv_h )& \le C \norm{\bu_{h}}_{E_h}  \norm{\bv_h}_{E_h},
\\ b_{h}(\bu_{h},\bv_h ) & \le C \norm{\bu_{h}}_{E_h}  \norm{\bv_h}_{E_h},
\end{align*}
for all $\bu_{h}$ and $\bv_h$ in $E_h$.
This implies that $a_h$, $T_h$ and $b_{h}$ are continuous.
The continuity of the mapping $\mathcal{P}$ follows from the continuity of the discrete forms
$a_h(\cdot,\cdot )$, $T_h(\cdot,\cdot )$ and $b_{h}(\cdot,\cdot)$.

Our goal now is to show that
\begin{equation}\label{eq:disc-exist:final}
	[\mathcal{P}(\bu_{h}^{n+1}),\bu_{h}^{n+1}]>0 
	\quad \text{for $\norm{\bu^{n+1}_h}_{E_h}=r>0$},
\end{equation}
for a sufficiently large $r$.  We observe that
\begin{equation}\label{eq:disc-exist:3}
\begin{split}
[\mathcal{P}(\bu^{n+1}_{h}),\bu^{n+1}_{h}]=& \frac{1}{\Delta t}\si \sum_{K \in \Omega_R }m(K)
\abs{u^{n+1}_{i,K}}^2
+a_h(\bu^{n+1}_{h},\bu^{n+1}_{h})
\\&\quad +b_{h}(\bu^{n+1}_{h},\bu^{n+1}_{h})
-\frac{1}{\Delta t}\si \sum_{K \in \Omega_R }m(K) u^{n}_{i,K}u^{n+1}_{i,K}.
\end{split}
\end{equation}
It follows that from the definition of $\sigma$ and \eqref{eq:disc-exist:3}
and Young's inequality that
\begin{equation*}
\begin{split}
&[\mathcal{P}(\bu^{n+1}_{h}),\bu^{n+1}_{h}]
\\ & \ge \frac{1}{\Delta t}
\si\sum_{K \in \Omega_R }m(K)\abs{u^{n+1}_{i,K}}^2
+\si M_i \norm{u^{n+1}_{i,h}}^2_{H_h(\Om)}\\
&\quad +\sum_{K \in \Omega_R }m(K) \Bigl(\sigma(u^{n+1,+}_{1,K},u^{n+1,+}_{2,K},u^{n+1,+}_{3,K})
u^{n+1,+}_{1,K}+\mu \abs{u^{n+1}_{1,K}}^2\Bigl)
\\
&\quad-\sum_{K \in \Omega_R }m(K) 
\Bigl(\sigma(u^{n,+}_{1,K},u^{n+1,+}_{2,K},u^{n+1,+}_{3,K})-(\gamma+\mu) u^{n+1}_{2,K}\Bigl)
u^{n+1}_{2,K},\\
&\quad-\gamma \sum_{K \in \Omega_R }m(K) u^{n}_{2,K}u^{n+1}_{3,K}
-\frac{1}{\Delta t}\si \sum_{K \in \Omega_R }m(K) u^{n}_{i,K}u^{n+1}_{i,K}
\\ & \ge\frac{1}{\Delta t}
\si\sum_{K \in \Omega_R }m(K)\abs{u^{n+1}_{i,K}}^2
+\si M_i \norm{u^{n+1}_{i,h}}^2_{H_h(\Om)}
\\
&\quad -\frac{1}{8 \Delta t}\sum_{K \in \Omega_R }m(K)\abs{u^{n+1}_{2,K}}^2
-C_1(\Delta t,\alpha)\sum_{K \in \Omega_R }m(K)\abs{u^{n,+}_{1,K}}^2\\
&\quad -\frac{1}{8 \Delta t}\sum_{K \in \Omega_R }m(K)\abs{u^{n+1}_{3,K}}^2
-C_2(\Delta t,\gamma)\sum_{K \in \Omega_R }m(K)\abs{u^{n}_{2,K}}^2\\
&\quad -\si \frac{1}{8 \Delta t}\sum_{K \in \Omega_R }m(K)\abs{u^{n+1}_{i,K}}^2
-C_3(\Delta t)\si\sum_{K \in \Omega_R }m(K)\abs{u^{n}_{i,K}}^2.
\end{split}
\end{equation*}
This implies that
\begin{equation}\label{eq:disc-exist:5}
\begin{split}
& [\mathcal{P}(u^{n+1}_{1,h}),u^{n+1}_{1,h}] 
\\ & \ge \min{\Bigl\{\frac{3}{4 \Delta t},M_1,M_2,M_3\Bigl\}}
\norm{\bu^{n+1}_{h}}^2_{E_h}
\\ & \qquad \qquad \qquad \qquad 
-2\max{\Bigl\{C_1(\Delta t,\alpha),C_2(\Delta t,\gamma),C_3(\Delta t)\Bigl\}}
\si\sum_{K \in \Omega_R }m(K)\abs{u^{n}_{i,K}}^2.
\end{split}
\end{equation}
Finally, for given $u^{n}_{1,h}$, $u^{n}_{2,h}$ and $u^{n}_{3,h}$,
we deduce from \eqref{eq:disc-exist:5} that \eqref{eq:disc-exist:final} holds for $r$ large enough
(recall that $\norm{\bu^{n+1}_{h}}_{E_h}=r$).
Hence, we obtain the existence of at least one solution to the scheme 
\eqref{prob:init}-\eqref{S3-discr}.

\end{proof}

\subsection{Nonnegativity}

We have the following lemma.
\begin{lem}
\label{lem:nonnegativity}
Let $(u_{1,K}^{n},u_{2,K}^{n},u_{3,K}^{n})_{K \in \Omega_R,n \in  \{0,\ldots,N\}}$
be a solution of the finite volume scheme \eqref{prob:init}, \eqref{S1-discr}, \eqref{S2-discr}
and \eqref{S3-discr}.
Then, $(u_{1,K}^{n},u_{2,K}^{n},u_{3,K}^{n})_{K \in \Omega_R,n \in  \{0,\ldots,N\}}$ is nonnegative.
\end{lem}

\begin{proof}
Multiplying \eqref{S1-discr} by $ -\Delta
t {u_{1,K}^{n+1}}^-$, we find that
\begin{equation}
   \label{nonnegat:I}
   \begin{split}
      &-m(K){u_{1,K}^{n+1}}^-(u^{n+1}_{1,K}-u^{n}_{1,K})
	 + a_1\Bigl(\sum_{K_0 \in \Om_h}u^{n}_{1,K_0}\Bigl) \Delta t \sum_{L \in N(K) }
     \frac{m(\sigma_{K,L})}{d(K,L)}(u^{n+1}_{1,L}-u^{n+1}_{1,K}){u_{1,K}^{n+1}}^-
     \\ &\qquad \qquad \qquad 
     -m(K)\Delta t  \Bigl(\sigma(u^{n+1,+}_{1,K},u^{n+1,+}_{2,K},u^{n+1,+}_{3,K}){u_{1,K}^{n+1}}^-
     +\mu u^{n+1}_{1,K}\Bigl) {u_{1,K}^{n+1}}^-=0.
   \end{split}
\end{equation}
We know that $u_{K}^{n+1}={u_{K}^{n+1}}^+-{u_{K}^{n+1}}^-$ and $(a^+-b^+)(a^--b^-)\leq 0$
for all $a,b\in\R$.
With this, we deduce 
\begin{equation}\label{nonnegat:I:1}\begin{split}
  &a_1\Bigl(\sum_{K_0 \in \Om_h}u^{n}_{1,K_0}\Bigl)\sn \Delta t \sum_{K \in \Omega_R} \sum_{L \in N(K) }
     \frac{m(\sigma_{K,L})}{d(K,L)}(u^{n+1}_{1,L}-u^{n+1}_{1,K}){u_{1,K}^{n+1}}^-\\
     &\qquad =-\frac{a_1\Bigl(\displaystyle\sum_{K_0 \in \Om_h}u^{n}_{1,K_0}\Bigl)}{2}
     \sn \Delta t \sum_{K \in \Omega_R} \sum_{L \in N(K) }
     (u^{n+1}_{1,L}-u^{n+1}_{1,K})({u_{1,L}^{n+1}}^--{u_{1,K}^{n+1}}^-)\\
     &\qquad \geq \frac{a_1\Bigl(\displaystyle\sum_{K_0 \in \Om_h}u^{n}_{1,K_0}\Bigl)}{2}
     \sn \Delta t \sum_{K \in \Omega_R} \sum_{L \in N(K) }
     \abs{{u_{1,L}^{n+1}}^--{u_{1,K}^{n+1}}^-}^2\geq 0,
\end{split}\end{equation}
and
\begin{equation}\begin{split}\label{nonnegat:I:2}
    & \sn  \Delta t  \sum_{K \in \Omega_R} m(K)
    \Bigl(\sigma(u^{n+1,+}_{1,K},u^{n+1,+}_{2,K},u^{n+1,+}_{3,K})
    +\mu u^{n+1}_{1,K} \Bigl){u_{1,K}^{n+1}}^-\\
    &\qquad \qquad \qquad \qquad \qquad \qquad 
    =-\mu \sn  \Delta t  \sum_{K \in \Omega_R} m(K)
\abs{{u_{1,K}^{n+1}}^-}^2\leq 0.
\end{split}\end{equation}
Let $f \in C^2$ function. By using a Taylor expansion we find
\begin{equation}\label{Taylor}
f(b)=f(a)+f'(a)(b-a)+\frac{1}{2}f''(\xi)(b-a)^2,
\end{equation}
for some $\xi$ between $a$ and $b$. Using the Taylor expansion
\eqref{Taylor} on the sequence $f(u_{1,K}^{n+1})$
with $\disp f(\rho)=\int_0^{\rho^-}s \ds$,
$a=u_{1,K}^{n+1}$ and $b=u_{1,K}^{n}$. We find
\begin{equation}\label{nonnegat:I:3}
{u_{1,K}^{n+1}}^-(u^{n+1}_{1,K}-u^{n}_{1,K})=\frac{\abs{{u_{1,K}^{n}}^-}^2}{2}
-\frac{\abs{{u_{1,K}^{n+1}}^-}^2}{2}
-\frac{1}{2}f''(\xi)\left(u^{n+1}_{1,K}-u^{n}_{1,K}\right)^2.
\end{equation}
We observe from the definition of $f$ that $f''(\rho)\geq 0$, which
implies
\begin{equation}\label{nonnegat:I:4}
{u_{1,K}^{n+1}}^-(u^{n+1}_{1,K}-u^{n}_{1,K})\le \frac{\abs{{u_{1,K}^{n}}^-}^2}{2}
-\frac{\abs{{u_{1,K}^{n+1}}^-}^2}{2}.
\end{equation}
Now, using \eqref{nonnegat:I:1}-\eqref{nonnegat:I:4} to deduce from
\eqref{nonnegat:I}
\begin{equation}
   \label{nonnegat-bis:I}
   \begin{split}
\sn \Bigl(\frac{\abs{{u_{1,K}^{n+1}}^-}^2}{2}
&-\frac{\abs{{u_{1,K}^{n}}^-}^2}{2}\Bigl)
+\frac{a_1\Bigl(\displaystyle\sum_{K_0 \in \Om_h}u^{n}_{1,K_0}\Bigl)}{2}
\sn \Delta t \abs{{u^{n+1}_{1,L}}^--{u^{n+1}_{1,K}}^-}^2\\
&\qquad \qquad \qquad \qquad 
+\mu \sn  \Delta t  \sum_{K \in \Omega_R} m(K) \abs{{u_{1,K}^{n+1}}^-}^2\le 0 .
   \end{split}
\end{equation}
This implies that
\begin{equation}
   \label{nonnegat-bis:I:1}
\frac{1}{2} 
\biggl(\abs{{u_{1,K}^{N}}^-}^2-\abs{{u_{1,K}^{0}}^-}^2\biggl)
\le 0 .
\end{equation}
Note that \eqref{nonnegat-bis:I:1} is also true if we replace $N$
by $n_0\in\{1,\dots,N\}$, so we have established
\begin{equation}
   \label{nonnegat-bis:I:1}
\abs{{u_{1,K}^{n_0}}^-}^2\le \abs{{u_{1,K}^{0}}^-}^2.
\end{equation}
Since $u_{1,K}^{0}$ is nonnegative, the result is
${u_{1,K}^{n+1}}^-$ for all $0\le n \le N-1$ and all $K \in \Omega_R$. 

On the other hand, multiplying \eqref{S2-discr} by $ -\Delta
t {u_{2,K}^{n+1}}^-$, and along 
 the
same lines as $u_{1,K}^{n+1}$, we obtain the nonnegativity
of discrete solutions $u_{2,K}^{n+1}$ for all $0\le n \le N-1$ and all $K \in \Omega_R$.
Finally, the nonnegativity of $u_{3,K}^{n+1}$, is given by a maximum principle proved in \cite{Eymard_etal_II:2000}.
In fact the reaction terms $m(K)\,\gamma\,u_{2,K}^{n}$ of the equation \eqref{S3-discr} do not depend on $u_{3,K}^{n+1}$
and it is a nonnegative term.  We assume the nonnegativity of 
$(u_{1,K}^{n},u_{2,K}^{n},u_{3,K}^{n})$ and we apply the discrete maximum principle for the
third equation \eqref{S3-discr}
in order to prove the  nonnegativity of $u_{3,K}^{n+1}$.
Then, using an induction on $n$ yields, we conclude.
\end{proof}

\section{A priori estimates}\label{sec:basic-apriori}

The goal is to establish several a priori (discrete energy) estimates 
for the finite volume scheme, which eventually will imply 
the desired convergence results.

\begin{prop}\label{prop:LPBV}
Let $(u_{i,K}^{n})_{K \in \Omega_R,n \in  \{0,\ldots,N\}}$, $i=1,2,3$,
be a solution of the finite volume scheme \eqref{prob:init}-\eqref{S3-discr}. 
Then there exist constants $C_{1},C_{2},C_{3}>0$, depending
on $\Omega$, $T$, $u_{i,0}$ and $\alpha$ such that
\begin{equation}\label{est:L2-norm}
\max_{n \in \{1,\dots,N\}}\sum_{K\in \Omega_{R}}m(K)\abs{u_{i,K}^{n}}^2
\le C_{1},
\end{equation}
\begin{equation}\label{est:grad-norm}
\frac{1}{2}\sn \Delta t \sum_{K\in \Omega_R}\sum_{L\in N(K)}
\frac{m(\sigma_{K,L})}{d(K,L)}
\abs{u_{i,K}^{n+1}-u_{i,L}^{n+1}}^2\le C_{2},
\end{equation}
and 
\begin{equation}\label{est:H1-norm}
\sn \Delta t \sum_{K\in \Omega_{R}}m(K)\Bigl(\abs{\sigma(u_{1,K}^{n},u_{2,K}^{n+1},u_{3,K}^{n+1})}^2
+\abs{\sigma(u_{1,K}^{n+1},u_{2,K}^{n+1},u_{3,K}^{n+1})}^2\Bigl)
\le C_{3}.
\end{equation}
for $i=1,2,3$.
\end{prop}

\begin{proof}
We multiply \eqref{S1-discr}, \eqref{S2-discr} and \eqref{S3-discr} by $\Delta t u^{n+1}_{1,K}$,
$\Delta t u^{n+1}_{2,K}$ and $\Delta t u^{n+1}_{3,K}$, respectively, and add together the outcomes.
Summing the resulting equation over $K$ and $n$ yields
$$
E_1 + E_2 + E_3 = E_4,
$$
where
$$
E_1 =\si \sn \sum_{K \in \Omega_R} m(K)(u^{n+1}_{i,K}-u^{n}_{i,K}) u^{n+1}_{i,K},
$$
$$
E_2 = -\si\sn \Delta t \,
a_i\Bigl(\sum_{K_0 \in \Om_h}u^{n}_{i,K_0}\Bigl)\sum_{K \in \Omega_R} \sum_{L \in N(K)}
\frac{m(\sigma_{K,L})}{d(K,L)}
(u^{n+1}_{i,L}-u^{n+1}_{i,K})u^{n+1}_{i,K},
$$
$$
E_3=\sn \Delta t \sum_{K \in \Omega_R}  m(K)
\Biggl(\Bigl(\sigma(u^{n+1}_{1,K},u^{n+1}_{2,K},u^{n+1}_{3,K})+\mu u^{n+1}_{1,K}\Bigl)u^{n+1}_{1,K}
+(\gamma+\mu) \abs{u^{n+1}_{2,K}}^2\Biggl),
$$
$$
E_4=\sn \Delta t \sum_{K \in \Omega_R} 
m(K)\Bigl(\sigma(u^{n}_{1,K},u^{n+1}_{2,K},u^{n+1}_{3,K})u^{n+1}_{2,K}\\
+\gamma u^{n}_{2,K}u^{n+1}_{3,K}\Bigl).
$$
From the inequality ``$a(a-b)\ge \frac{1}{2}(a^2-b^2)$", we obtain
\begin{equation*}
\begin{split}
E_1 &= \si \sn \sum_{K \in \Omega_R} m(K)(u^{n+1}_{i,K}-u^{n}_{i,K}) u^{n+1}_{i,K}\\
&\ge \frac{1}{2}\si \sn \sum_{K \in \Omega_R} m(K)
\left(\abs{u^{n+1}_{i,K}}^2-\abs{u^{n}_{i,K}}^2\right)  \\
&= \frac{1}{2}\si \sum_{K \in \Omega_R} m(K)\left(\abs{u^{N}_{i,K}}^2-\abs{u^{0}_{i,K}}^2\right) .
\end{split}
\end{equation*}
Gathering by edges, we obtain
\begin{equation*}
      \begin{split}
E_2 &=-\si \sn \Delta t \, a_i\Bigl(\sum_{K_0 \in \Om_h}u^{n}_{i,K_0}\Bigl)
\sum_{K \in \Omega_R} \sum_{L \in N(K)}
\frac{m(\sigma_{K,L})}{d(K,L)}
(u^{n+1}_{i,L}-u^{n+1}_{i,K})u^{n+1}_{i,K}\\
&\ge \si \frac{M_i}{2}\sn \Delta t \sum_{K \in \Omega_R} \sum_{L \in N(K)}
\frac{m(\sigma_{K,L})}{d(K,L)}
\abs{u^{n+1}_{i,K}-u^{n+1}_{i,L}}^2.
      \end{split}
\end{equation*}
Observe that from nonnegativity of
$(u_{i,K}^{n})_{K \in \Omega_R,n \in  \{0,\ldots,N\}}$ for $i=1,2,3$, we get
\begin{equation*}
E_3\ge 0. 
\end{equation*}
We use Young's inequality to deduce
\begin{equation*}
  \begin{split}
E_4=&\sn \Delta t \sum_{K \in \Omega_R} 
m(K)\Bigl(\sigma(u^{n}_{1,K},u^{n+1}_{2,K},u^{n+1}_{3,K})u^{n+1}_{2,K}
+\gamma u^{n}_{2,K}u^{n+1}_{3,K}\Bigl)\\
\le &\alpha \sn \Delta t \sum_{K \in \Omega_R}  m(K) \abs{u^{n+1}_{2,K}}^2
+\frac{\gamma}{2} \sn \Delta t \sum_{K \in \Omega_R}  m(K) \abs{u^{n}_{2,K}}^2
\\
&\qquad \qquad \qquad \qquad 
+\frac{\gamma}{2} \sn \Delta t \sum_{K \in \Omega_R}  m(K) \abs{u^{n+1}_{3,K}}^2\\
\le &(\alpha+\frac{\gamma}{2}) \sn \Delta t \sum_{K \in \Omega_R}  m(K) \abs{u^{n+1}_{2,K}}^2
+\frac{\gamma}{2}\Delta t \sum_{K \in \Omega_R}  m(K) \abs{u^{0}_{2,K}}^2\\
&\qquad \qquad \qquad \qquad +\frac{\gamma}{2} \sn \Delta t \sum_{K \in \Omega_R}  m(K) \abs{u^{n+1}_{3,K}}^2.
 \end{split}
\end{equation*} 
Collecting the previous inequalities we obtain
\begin{equation}\label{eq:collect:1}
      \begin{split}
&\frac{1}{2}\si \sum_{K \in \Omega_R} m(K)(\abs{u^{N}_{i,K}}^2-\abs{u^{0}_{i,K}}^2)
\\
&\quad +\si \frac{M_i}{2}
\sn \Delta t \sum_{K \in \Omega_R} \sum_{L \in N(K)}
\frac{m(\sigma_{K,L})}{d(K,L)}
\abs{u^{n+1}_{i,K}-u^{n+1}_{i,L})}^2\\
\le& (\alpha+\frac{\gamma}{2}) \sn \Delta t \sum_{K \in \Omega_R}  m(K) \abs{u^{n+1}_{2,K}}^2
+\frac{\gamma}{2}\Delta t \sum_{K \in \Omega_R}  m(K) \abs{u^{0}_{2,K}}^2\\
&\qquad \qquad \qquad \qquad +\frac{\gamma}{2} \sn \Delta t \sum_{K \in \Omega_R}  m(K) \abs{u^{n+1}_{3,K}}^2,
      \end{split}
\end{equation}      
which implies
\begin{equation}\label{eq:collect:2}
      \begin{split}
&\si \sum_{K \in \Omega_R} m(K)\abs{u^{N}_{i,K}}^2\\
&\qquad\qquad \le 
\si \sum_{K \in \Omega_R} m(K)\abs{u^{0}_{i,K}}^2
+(\alpha+\frac{\gamma}{2}) \sn \Delta t \sum_{K \in \Omega_R}  m(K) \abs{u^{n+1}_{2,K}}^2
\\
&\qquad\qquad \qquad +\frac{\gamma}{2}\Delta t \sum_{K \in \Omega_R}  m(K) \abs{u^{0}_{2,K}}^2
+\frac{\gamma}{2} \sn \Delta t \sum_{K \in \Omega_R}  m(K) \abs{u^{n+1}_{3,K}}^2.
      \end{split}
\end{equation}      
Clearly,  
\begin{equation*}
\sum_{K \in \Omega_R} m(K)\abs{u^{0}_{i,K}}^2
\le \norm{u_{i,0}}^2_{L^2(\Om)} \text{ for }i=1,2,3.
\end{equation*}     
In view of \eqref{eq:collect:2}, this implies that there exist constants
$C_4,C_5>0$ such that
\begin{equation}\label{eq:collect:3-bis}
\si\sum_{K \in \Omega_R} m(K)\abs{u^{N}_{i,K}}^2\le
C_4+ C_5\si \sn \Delta t \sum_{K \in \Omega_R} 
m(K) \abs{u^{n+1}_{i,K}}^2.
\end{equation}      

Note that \eqref{eq:collect:3-bis} is also true if we replace $N$
by $n_0\in\{1,\dots,N\}$, so we have established
\begin{equation}\label{eq:collect:3}
\si\sum_{K \in \Omega_R} m(K)\abs{u^{n_0}_{i,K}}^2\le
C_4+ C_5\si \sum_{n=0}^{n_0-1} \Delta t \sum_{K \in \Omega_R} 
m(K) \abs{u^{n+1}_{i,K}}^2.
\end{equation}     
By the discrete Gronwall inequality (see e.g. \cite{Fairweather}), we obtain from \eqref{eq:collect:3}
\begin{equation}\label{eq:collect:4}
\si \sum_{K \in \Omega_R} m(K)\abs{u^{n_0}_{i,K}}^2\le
C_6,
\end{equation}      
for any $n_0\in\{1,\dots,N\}$ and some constant $C_6>0$. Then
\begin{equation*}
\max_{n \in \{1,\dots,N\}}\si \sum_{K \in \Omega_R} m(K)\abs{u^{n}_{i,K}}^2\le
C_6.
\end{equation*}      
Moreover, we obtain from \eqref{eq:collect:1} and \eqref{eq:collect:4} the existence of 
a constant $C_7>0$ such that
\begin{equation*}
\si \sn \Delta t \sum_{K \in \Omega_R} \sum_{L \in N(K)}
\frac{m(\sigma_{K,L})}{d(K,L)}
\abs{u^{n+1}_{i,K}-u^{n+1}_{i,L}}^2\le C_7.
\end{equation*} 
Finally a consequence of \eqref{def-sigma-1} and \eqref{est:L2-norm}
is that
\begin{equation*}
\norm{\sigma(u_{1,h},u_{3,h},u_{3,h})}_{L^2(Q_{T})} 
\le C_{8},
\end{equation*}
for some constant $C_8>0$. 

\end{proof}

\section{Space and time translation estimates}\label{sec-compactness}

In this section we derive estimates on differences of space and time translates of the
function $v_{h}$ which imply that the sequence
$v_{h}$ is relatively compact in $L^2(Q_T)$.

\begin{lem}\label{Space-Time-translate}
There exists a constant $C>0$ depending on
$\Omega$, $T$, $u_{1,0}$, $u_{2,0}$, $u_{3,0}$, $\alpha$, $\gamma$ and $\mu$ such that 
\begin{equation}\label{Space-translate}
\iint_{\Om' \times (0,T)}\abs{u_{i,h}(t,x+y)-u_{i,h}(t,x)}^2\dx \dt
\le C\abs{y}(\abs{y}+2h),\quad i=1,2,3,
\end{equation}
for all $y \in \R^3$ with $\Om'=\{x \in \Om, \, [x,x+y]\subset \Om\}$, and
\begin{equation}\label{Time-translate}
\iint_{\Om \times
(0,T-\tau)}\abs{u_{i,h}(t+\tau,x)-u_{i,h}(t,x)}^2\dx \dt \le C(\tau+\Delta t),\quad i=1,2,3,
\end{equation}
for all $\tau\in (0,T)$.
\end{lem}

\begin{proof}
The proof is similar to that found in, e.g, \cite{Eymard_etal_II:2000}.

{\it Proof of \eqref{Space-translate}}.
Let $y \in \R^3$, $x \in \Om'$, and $L \in N(K)$. We set 
$$
\chi_{\sigma_{K,L}}=
\begin{cases}
	1, & \text{if the line segment $[x,x+y]$ intersects $\sigma_{K,L}$, $K$ and $L$},\\
	0, & \text{otherwise}. 
\end{cases}
$$
Next, the value $c_{\sigma_{K,L}}$ is defined by
$\displaystyle c_{\sigma_{K,L}}=\frac{y}{\abs{y}}\cdot \eta_{K,L}$
with $c_{\sigma_{K,L}}>0$.
Observe that
\begin{equation}\label{est:space}
\int_{\Om'}\chi_{\sigma_{K,L}}(x) \dx \le m(\sigma_{K,L}) \abs{y} c_{\sigma_{K,L}}.
\end{equation}
Using this, we obtain 
\begin{equation*}
      \begin{split}
    \abs{u_{1,h}(t,x+y)-u_{1,h}(t,x)}\le
   \sum_{\sigma_{K,L}} \chi_{\sigma_{K,L}}(x)
\abs{u_{1,L}^{n+1}-u_{1,K}^{n+1}}.
      \end{split}
\end{equation*}
To simplify the notation, we write 
$$
\sum_{\sigma_{K,L}}\quad \text{instead of} \quad 
\sum_{\{(K,L)\in \Omega^2_R,\,K\ne L,\,m(\sigma_{K,L})\ne 0\}}.
$$
By the Cauchy-Schwarz inequality, we get
\begin{equation}\label{est:space:2}
\begin{split}
\abs{u_{1,h}(t,x+y)-u_{1,h}(t,x)}^2 
\le & \sum_{\sigma_{K,L}} \chi_{\sigma_{K,L}}(x)c_{\sigma_{K,L}}
d(K,L) \\ &\quad \times 
\sum_{\sigma_{K,L}} \frac{\abs{u_{1,L}^{n+1}-u_{1,K}^{n+1}}^2}
{c_{\sigma_{K,L}}d(K,L)} \chi_{\sigma_{K,L}}(x).
\end{split}
\end{equation}
Note that
\begin{equation}\label{est:space:3}
\sum_{\sigma_{K,L}} \chi_{\sigma_{K,L}}(x)
c_{\sigma_{K,L}}d(K,L)\le \abs{y}+2h.
\end{equation}
Using \eqref{est:space}, \eqref{est:space:2}, and \eqref{est:space:3}, we deduce
\begin{equation}\label{est:space:4}
      \begin{split}
&\iint_{(0,T)\times\Om'}\abs{u_{1,h}(t,x+y)-u_{1,h}(t,x)}^2\dx \\
&\qquad \le  (\abs{y}+2h)\sn \Delta t
\sum_{\sigma_{K,L}}\frac{\abs{u_{1,L}^{n+1}-u_{1,K}^{n+1}}^2}
{c_{\sigma_{K,L}}d(K,L)}\int_{\Om'} \chi_{\sigma_{K,L}}(x)\dx\\
&\qquad   \le \abs{y}(\abs{y}+2h)\sn \Delta t \sum_{\sigma_{K,L}}
\frac{m(\sigma_{K,L})}{d(K,L)}\abs{u_{1,L}^{n+1}-u_{1,K}^{n+1}}^2.
      \end{split}
\end{equation}
Then, from \eqref{est:grad-norm} and \eqref{est:space:4}, we 
deduce \eqref{Space-translate}. 

{\it Proof of \eqref{Time-translate}}. 
Let $\tau \in (0,T)$ and $t \in (0,T-\tau)$. We have
\begin{equation*}
      \begin{split}
\mathcal{B}(t)=&\int_{\Om} \abs{u_{1,h}(t+\tau,x)-u_{1,h}(t,x)}^2\dx.
      \end{split}
\end{equation*}
Set $n_0(t)=E(t/\Delta t)$ and $n_1(t)=E((t+\tau)/\Delta t)$, where 
\begin{equation*}
\text{$E(x)=n$ for $x\in [n, n+1)$, $n \in \N$.}
\end{equation*}
We get
\begin{equation*}
  \mathcal{B}(t)= \sum_{ K \in \Omega_R}m(K)
  \abs{u_{1,K}^{n_1(t)}-u_{1,K}^{n_0(t)}}^2.
\end{equation*}
This implies
\begin{equation*}\label{Time-translate:3}
  \mathcal{B}(t)= \sum_{ L \in \Om_R}m(K)
 \Bigl (u_{1,K}^{n_1(t)}-u_{1,K}^{n_0(t)}\Bigl)
 \sum_{{t<(n+1) \Delta t<t+\tau}}
m(K)(u_{1,K}^{n+1}-u_{1,K}^{n}).
 \end{equation*}
Using the scheme \eqref{S1-discr}, we obtain
\begin{equation}\label{Time-translate:4}
      \begin{split}
\mathcal{B}(t)=& \sum_{{t<(n+1)\Delta t<t+\tau}} \Delta t
\sum_{ K \in \Omega_R}
  \Bigl(u_{1,K}^{n_1(t)}-u_{1,K}^{n_0(t)}\Bigl)\\
&\qquad \qquad
\times\Biggl(a_1\Bigl(\sum_{K_0 \in \Om_h}u^{n}_{1,K_0}\Bigl)
\sum_{L \in N(K) }
m(\sigma_{K,L}) \frac{u_{1,L}^{n+1}-u_{1,K}^{n+1}}{d(K,L)}\\
&  \qquad \qquad  \qquad \qquad
-m(K)\Bigl(\sigma(u^{n+1}_{1,K},u^{n+1}_{2,K},u^{n+1}_{3,K})
+\mu u^{n+1}_{1,K}\Bigl)\Biggl).
     \end{split}
\end{equation}

We observe that we can rewrite \eqref{Time-translate:4} as
\begin{equation*}
      \begin{split}
 \mathcal{B}(t)=& \frac{1}{2 }\sum_{{t< (n+1)\Delta t<t+\tau}} \Delta t\,
a_1\Bigl(\sum_{K_0 \in \Om_h}u^{n}_{1,K_0}\Bigl)
\sum_{ K \in \Omega_R}\sum_{L \in N(K) }\frac{m(\sigma_{K,L})}{d(K,L)}\\
&\times \Biggl[ \Bigl(u_{1,K}^{n+1}-u_{1,L}^{n+1}\Bigl)
\Bigl(u_{1,K}^{n_1(t)}-u_{1,L}^{n_1(t)}\Bigl)
+\Bigl(u_{1,L}^{n+1}-u_{1,K}^{n+1}\Bigl)
\Bigl(u_{1,L}^{n_0(t)}-u_{1,K}^{n_0(t)}\Bigl)\Biggl]\\
& \qquad \quad +\sum_{{t<(n+1)\Delta t<t+\tau}} \Delta t
\sum_{ K \in \Omega_R}
  \Bigl(u_{1,K}^{n_1(t)}-u_{1,K}^{n_0(t)}\Bigl)\\
&\qquad \qquad \qquad
\times \Biggl(-m(K)\Bigl(\sigma(u^{n+1}_{1,K},u^{n+1}_{2,K},u^{n+1}_{3,K})
+\mu u^{n+1}_{1,K}\Bigl)\Biggl).
     \end{split}
\end{equation*}
We use the basic inequality ``$ab \le \frac{1}{2}a^2+\frac{1}{2}b^2$" to deduce
\begin{equation*}
\mathcal{B}(t)\le \frac{1}{2}
\Bigl(\mathcal{B}_1(t)+\frac{1}{2}\mathcal{B}_2(t)
+\frac{1}{2}\mathcal{B}_3(t)\Bigl)
+\mathcal{B}_4(t),
\end{equation*}
with
\begin{equation*}\begin{split}
\mathcal{B}_1(t)&=\sum_{{t< (n+1) \Delta t<t+\tau}} \Delta t\;
a_1\Bigl(\sum_{K_0 \in \Om_h}u^{n}_{1,K_0}\Bigl)
\sum_{ K \in \Omega_R}\sum_{L \in N(K) }\frac{m(\sigma_{K,L})}{d(K,L)}
\abs{u_{1,K}^{n+1}-u_{1,L}^{n+1}}^2,\\
\mathcal{B}_2(t)&=\sum_{{t< (n+1)\Delta t<t+\tau}} \Delta t\;
a_1\Bigl(\sum_{K_0 \in \Om_h}u^{n}_{1,K_0}\Bigl)
\sum_{ K \in \Omega_R}\sum_{L \in N(K) }\frac{m(\sigma_{K,L})}{d(K,L)}
\abs{u_{1,K}^{n_1(t)}-u_{1,L}^{n_1(t)}}^2,\\
\mathcal{B}_3(t)&=\sum_{{t<(n+1) \Delta t<t+\tau}} \Delta t\;
a_1\Bigl(\sum_{K_0 \in \Om_h}u^{n}_{1,K_0}\Bigl)
\sum_{ K \in \Omega_R}\sum_{L \in N(K) }\frac{m(\sigma_{K,L})}{d(K,L)}
\abs{u_{1,K}^{n_0(t)}-u_{1,L}^{n_0(t)}}^2,\\
\mathcal{B}_4(t)&=-\sum_{{t<(n+1)\Delta t<t+\tau}} \Delta t
\sum_{ K \in \Omega_R}m(K)
  (u_{1,K}^{n_1(t)}-u_{1,K}^{n_0(t)})\\
&\qquad \qquad \qquad \qquad \qquad \qquad 
\times \Bigl(\sigma(u^{n+1}_{1,K},u^{n+1}_{2,K},u^{n+1}_{3,K})
+\mu u^{n+1}_{1,K}\Bigl).
\end{split}\end{equation*}

Next, we introduce the characteristic function $\chi$ defined by
$\chi(n,t_1,t_2)= 1$ if $t_1 < (n+1) \Delta t \le t_2$ and $\chi(n, t_1, t_2) = 0$ otherwise. 
Then we have for any sequence $(a^n)_{n \in \N}$ of nonnegative numbers that 
\begin{equation}\label{lem:time:1}
\int_0^{T-\tau}\sum_{{t< (n+1) \Delta t<t+\tau}} a^n \dt \le
\sum_{n=0}^{E(\frac{T}{\Delta t})} a^n \int_0^{T-\tau}\chi(n, t, t+\tau)\dt \le
\tau \sum_{n=0}^{E(\frac{T}{\Delta t})} a^n
\end{equation}
and for any $\xi \in [0,\tau]$
\begin{equation}\label{lem:time:2}
\begin{split}
&\int_0^{T-\tau}\sum_{{t< (n+1) \Delta t<t+\tau}} a^{E((t+\xi)/\Delta t)}\\
& \le  
\sum_{m=0}^{E(\frac{T}{\Delta t})}\int_{m\Delta t}^{(m+1)\Delta t} a^m\sum_{n=0}^{E(\frac{T}{\Delta t})}
\chi(n, t-\xi, t-\xi+\tau)\dt\\
& =\sum_{m=0}^{E(\frac{T}{\Delta t})}\int_{0}^{\Delta t} a^m\sum_{n=0}^{E(\frac{T}{\Delta t})}
\chi(n-m, t-\xi, t-\xi+\tau)\dt\\
& =\sum_{m=0}^{E(\frac{T}{\Delta t})}\sum_{n=0}^{E(\frac{T}{\Delta t})}\int_{-n\Delta t}^{\Delta t-n\Delta t} a^m
\chi(-m, t-\xi, t-\xi+\tau)\dt\\
&  =\sum_{m=0}^{E(\frac{T}{\Delta t})}a^m \int_{\R}
\chi(-m, t-\xi, t-\xi+\tau)\dt\\
& =\sum_{m=0}^{E(\frac{T}{\Delta t})}a^m
\int_{\xi+\tau-(n+1)\Delta t}^{\xi-(n+1)\Delta t}\dt 
 = \tau  \sum_{n=0}^{E(\frac{T}{\Delta t})} a^n.
\end{split}
\end{equation}

From \eqref{lem:time:1}, we deduce
\begin{equation*}
     \begin{split}
\int^{T-\tau}_0\mathcal{B}_1(t)\dt &\le \sn \Delta t \int_0^{T-\tau}\chi(n, t, t+\tau)
\, a_1\Bigl(\sum_{K_0 \in \Om_h}u^{n}_{1,K_0}\Bigl)\\
&\qquad \times \sum_{ K \in \Omega_R}\sum_{L \in N(K) }\frac{m(\sigma_{K,L})}{d(K,L)}
\abs{u_{1,K}^{n+1}-u_{1,L}^{n+1}}^{2} \dt\\
&\le
\tau \sn \Delta t \, a_1\Bigl(\sum_{K_0 \in \Om_h}u^{n}_{1,K_0}\Bigl)
\sum_{(K,L)\in \Omega_R^2}\frac{m(\sigma_{K,L})}{d(K,L)}
\abs{u_{1,K}^{n+1}-u_{1,L}^{n+1}}^{2}.
     \end{split}
\end{equation*}
Using \eqref{est:grad-norm}, this implies that there exists a 
constant $C_{10}>0$ such that
\begin{equation*}
\int^{T-\tau}_0\mathcal{B}_1(t)\dt \le \tau C_{10}.
\end{equation*}
Now we consider $\mathcal{B}_2(t)$ and $\mathcal{B}_3(t)$.
We use \eqref{lem:time:2} with $\xi=\tau$ for $\mathcal{B}_2(t)$ and \eqref{lem:time:1}
for $\mathcal{B}_3(t)$ to obtain
$$
\int^{T-\tau}_0\mathcal{B}_2(t)\dt \le \tau \sn \Delta t \;
a_1\Bigl(\sum_{K_0 \in \Om_h}u^{n}_{1,K_0}\Bigl)
\sum_{ K \in \Omega_R}\sum_{L \in N(K) } \frac{m(\sigma_{K,L})}{d(K,L)}
\abs{u_{1,K}^{n+1}-u_{1,L}^{n+1}}^{2}
$$
and
$$
\int^{T-\tau}_0\mathcal{B}_3(t)\dt \le
\tau \sn \Delta t \;a_1\Bigl(\sum_{K_0 \in \Om_h}u^{n}_{1,K_0}\Bigl)
\sum_{ K \in \Omega_R}\sum_{L \in N(K) }\frac{m(\sigma_{K,L})}{d(K,L)}
\abs{u_{1,K}^{n+1}-u_{1,L}^{n+1}}^{2}.
$$
In view of \eqref{est:grad-norm}, we deduce
\begin{equation*}
  \int^{T-\tau}_0\mathcal{B}_2(t)\dt \le \tau C_{11}, \qquad
  \int^{T-\tau}_0\mathcal{B}_3(t)\dt \le \tau C_{12},
\end{equation*}
for some constants $C_{11},C_{12}>0$.

Next, we intergrate $\mathcal{B}_4(t)$ from $0$ to $T-\tau$. 
Using the inequalities ``$ab\le \frac{a^2}{2}+\frac{b^2}{2}$" and
``$(a-b)^2\le 2a^2+2b^2$", we obtain
\begin{equation*}
\begin{split}
 \mathcal{B}_4(t)&=   -\sum_{{t<(n+1)\Delta t<t+\tau}} \Delta t
    \sum_{ K \in \Omega_R}m(K)
    (u_{1,K}^{n_1(t)}-u_{1,K}^{n_0(t)})\Bigl(\sigma(u^{n+1}_{1,K},u^{n+1}_{2,K},u^{n+1}_{3,K})
+\mu u^{n+1}_{1,K}\Bigl)\\
    &\le \sum_{{t<(n+1)\Delta t<t+\tau}} \Delta t
    \sum_{ K \in \Omega_R}m(K)
    \Biggl(\abs{u_{1,K}^{n_1(t)}}^2+ \abs{u_{1,K}^{n_0(t)}}^2+
    C(\alpha,\mu)\frac{\abs{u_{1,K}^{n+1})}^2}{2}\Biggl).
  \end{split}
\end{equation*}
From \eqref{lem:time:1} and \eqref{est:L2-norm}, we have
$$
\int_0^{T-\tau}\sum_{{t<(n+1)\Delta t<t+\tau}} \Delta t
    \sum_{ K \in \Omega_R}m(K)
    \Biggl(\abs{v_K^{n_1(t)}}^2+ \abs{v_K^{n_0(t)}}^2\Biggl) \dt \le \tau C_{13},
$$
for some constant $C_{13}>0$. We use \eqref{est:L2-norm} to deduce
$$
\int_0^{T-\tau}\sum_{{t<(n+1)\Delta t<t+\tau}} \Delta t 
\sum_{ K \in \Omega_R}m(K) \frac{\abs{u_{1,K}^{n+1}}^2}{2}\le \tau C_{14},
$$
for some constant $C_{14}>0$. Hence
$$
\int_0^{T-\tau} \mathcal{B}_4(t)\le \tau C_{15},
$$
for some constant $C_{15}>0$.
Reasoning along the same lines for $u_{1,h}$ yield
\eqref{Space-translate} and \eqref{Time-translate} for $u_{2,h}$ and $u_{3,h}$.
This concludes the proof of the lemma.
\end{proof}

\section{Convergence of the finite volume scheme}\label{sec:conv}

The next lemma is a consequence of Lemma \ref{Space-Time-translate} 
and Kolmogorov's compactness criterion (see, e.g., \cite{Brezis}, Theorem IV.25).
\begin{lem}\label{lem-conv:1}
There exists a subsequence of $\bu_h=(u_{1,h},u_{2,h},u_{3,h})$, not relabeled, such 
that, as $h\to 0$,
\begin{equation}\label{limit-strong}
\begin{split}
& \text{(i)\; $u_{i,h}\to u_i$ strongly in $L^2(Q_T)$ and a.e. in $Q_T$},\\
& \text{(ii)\; $\Grad_h u_{i,h}\to \Grad u_i$ weakly in $(L^2(Q_T))^3$},
\\ 
&\text{(iii)\; $\sigma(u_{1,h},u_{2,h},u_{3,h})\to \sigma(u_1,u_2,u_3)$ strongly
in $L^1(Q_T)$},
\end{split}
\end{equation}
for $i=1,2,3$.
\end{lem}

\begin{proof}
The claim (i) in \eqref{limit-strong} follows from Lemma \ref{Space-Time-translate} 
and Kolmogorov's compactness criterion (see, e.g., \cite{Brezis}, Theorem IV.25).
The proof of the claim (ii) will be omitted since it is similar to that of Lemma 4.4 in
\cite{Chainais-1}, we refer to the proof of this lemma for more details. The claim (iii)
follows from Vitali theorem.
\end{proof}

Our final goal is to prove that the limit functions $u_1,u_2,u_3$
constructed in Lemma \ref{lem-conv:1} 
constitute a weak solution of the nonlocal system \eqref{S1}-\eqref{S3}. 

We start by verifying \eqref{eq1:def1}. 
Let $T$ be a fixed positive constant and
$\varphi_1\in \D([0,T)\times \overline{\Om})$.
We multiply the discrete equation \eqref{S1-discr} by
$\Delta t\varphi_1(t^n,x_K)$ for all $K \in \Omega_R$ and $n \in  \{0,\ldots,N\}$.
Summing the result over $K$ and $n$ yields
\begin{equation*}
T_1+T_2+T_3=0,
\end{equation*}
where
\begin{equation*}\begin{split}
T_1&=\sn \sum_{K\in \Omega_R}m(K)(u^{n+1}_{1,K}-u^{n}_{1,K})\varphi_1(t^n,x_K),\\
T_2&=-\sn \Delta t \,a_1\Bigl(\sum_{K_0 \in \Om_h}u^{n}_{1,K_0}\Bigl)
\sum_{K\in \Omega_R} \sum_{L\in N(K)}
\frac{m(\sigma_{K,L})}{d(K,L)}
(u^{n+1}_{1,L}-u^{n+1}_{1,K})\varphi_1(t^n,x_K),\\
T_3&=\sn \Delta t \sum_{K\in \Omega_R}m(K)
\Bigl(\sigma(u^{n+1}_{1,K},u^{n+1}_{2,K},u^{n+1}_{3,K})
+\mu u^{n+1}_{1,K}\Bigl)\varphi_1(t^n,x_K).
\end{split}\end{equation*}

Doing integration-by-parts, keeping in mind 
that $\varphi_1(T,x_K)=0$ for all $K \in \Omega_R$, we obtain
\begin{equation*}
\begin{split}
T_1&=-\sn \sum_{K\in \Omega_R}m(K)u^{n+1}_{1,K}(\varphi_1(t^{n+1},x_K)-\varphi_1(t^n,x_K))
-\sum_{K\in \Omega_R}m(K) u^{0}_{1,K}\varphi_1(0,x_K)\\
&=-\sn \sum_{K\in \Omega_R} \int_{t^n}^{t^{n+1}}\int_K u^{n+1}_{1,K}
\pt \varphi_1(t,x_K)\dx \dt-\sum_{K\in \Omega_R}\int_K u_{1,0}(x)\varphi_1(0,x_K)\dx\\
&=:-T_{1,1}-T_{1,2}.
\end{split}
\end{equation*}
Let us also introduce 
\begin{equation*}
\begin{split}
T_1^* &=-\sn \sum_{K\in \Omega_R} \int_{t^n}^{t^{n+1}}\int_K u^{n+1}_{1,K}
\pt \varphi_1(t,x)\dx \dt-\int_\Om u_{1,0}(x)\varphi_1(0,x)\dx\\
&=:-T_{1,1}^*-T_{1,2}^*.
\end{split}
\end{equation*}
Then
\begin{equation*}
\begin{split}
T_{1,2}-T_{1,2}^*&=\sum_{K\in \Omega_R}\int_K u_{1,0}(x)(\varphi_1(0,x_K)-\varphi_1(0,x))\dx.
\end{split}
\end{equation*}
From the regularity of $\varphi_i$, there exists a positive constant $C$ such that
$$
\abs{\varphi_1(0,x_K)-\varphi_1(0,x)}\le C\,h,
$$
which implies
\begin{equation}\label{conv:est:1}
\begin{split}
\abs{T_{1,2}-T_{1,2}^*}&\le C \, h\sum_{K\in \Omega_R}\int_K \abs{u_{1,0}(x)}\dx.
\end{split}
\end{equation}
Sending $h\to 0$ in \eqref{conv:est:1}, we get 
$$
\lim_{h \to 0}\abs{T_{1,2}-T_{1,2}^*}=0.
$$
Observe that
\begin{equation*}
\begin{split}
&T_{1,1}-T_{1,1}^*\\
& =\sn \sum_{K\in \Omega_R} \Biggl(\int_{t^n}^{t^{n+1}}\int_K u^{n+1}_{1,K}
\pt \varphi_1(t,x_K)\dx \dt-\int_{t^n}^{t^{n+1}}\int_K u^{n+1}_{1,K} 
\pt \varphi_1(t,x)\dx \dt\Biggl)\\
&=\sn \sum_{K\in \Omega_R} u^{n+1}_{1,K} \int_{t^n}^{t^{n+1}}\int_K
\Bigl( \pt \varphi_1(t,x_K)-\pt \varphi_1(t,x)\Bigl)\dx\dt.
\end{split}
\end{equation*}
Using the regularity of $\pt \varphi_1$ and H\"{o}lder's inequality, we obtain
\begin{equation*}
\begin{split}
\abs{T_{1,1}-T_{1,1}^*}&\le C(h)\Biggl(\sn \Delta t\sum_{K\in \Omega_R}
m(K)\abs{u^{n+1}_{1,K}}^2\Biggl)^{1/2},
\end{split}
\end{equation*}
where $C(h)>0$ is a function satisfying $C(h)\to 0$ as $h \to 0$.
From \eqref{prop:LPBV} we deduce
$$
\lim_{h \to 0}\abs{T_{1,1}-T_{1,1}^*}=0.
$$
Next, we define $I_D$ and $T_2^*$ by
\begin{equation*}
\begin{split}
&I_D= \int_0^T a_1\Bigl(\int_\Om u_1 \dx\Bigl)
\int_\Om \Grad u_1 \cdot \Grad \varphi_1\dx \dt,\\
&T_2^*=\int_0^T a_1\Bigl(\int_\Om u_{1,h} \dx\Bigl)
\int_\Om u_{1,h} \Delta \varphi_1\dx \dt.
\end{split}
\end{equation*}
Integration-by-parts yields
$$
I_D=-\int_0^T a_1\Bigl(\int_\Om u_1 \dx\Bigl)
\int_\Om u_1 \Delta \varphi_1\dx \dt.
$$
On the other hand,
using the convergence results of Lemma \ref{lem-conv:1}, and taking into account
the assumption \eqref{ass-ai}, it is easy
to prove that 
$a_1\left(\int_\Om u_{1,h} \dx\right) \to a_1\left(\int_\Om u_1 \dx\right) $
strongly in $L^2(0,T)$ and
$\int_\Om u_{1,h} \Delta \varphi_1\dx \to -\int_\Om \Grad u_1 \cdot \Grad \varphi_1\dx $
weakly in $L^2(0,T)$, as $h\to 0$.
Thus, there holds
$$
T_2^* \to -I_D \quad \text{as $h\to 0$}.
$$
Note that
\begin{equation*}
  \begin{split}
    T^*_2&=\sn a_1\Bigl(\sum_{K_0 \in \Om_h}u^{n}_{1,K_0}\Bigl)
\sum_{K\in \Omega_R} \sum_{L\in N(K)}
     u^{n+1}_{1,K}\int_{t_n}^{t_{n+1}}
\int_{\sigma_{K,L}}\Grad \varphi_1\cdot\eta_{K,L}d\gamma\\
&=-\frac{1}{2} \sn a_1\Bigl(\sum_{K_0 \in \Om_h}u^{n}_{1,K_0}\Bigl)
\sum_{K\in \Omega_R} \sum_{L\in N(K)}
 (u^{n+1}_{1,L}-u^{n+1}_{1,K})\int_{t_n}^{t_{n+1}}
\int_{\sigma_{K,L}}\Grad \varphi_1\cdot\eta_{K,L}d\gamma
  \end{split}
\end{equation*}
and
\begin{equation*}
  \begin{split}
    T_2&=-\sn \Delta t a_1\Bigl(\sum_{K_0 \in \Om_h}u^{n}_{1,K_0}\Bigl)
\sum_{K\in \Omega_R} \sum_{L\in N(K)}
   m(\sigma_{K,L})
\frac{u^{n+1}_{1,L}-u^{n+1}_{1,K}}{d(K,L)}\varphi_1(t^n,x_K)\\
&=\frac{1}{2} \sn a_1\Bigl(\sum_{K_0 \in \Om_h}u^{n}_{1,K_0}\Bigl)
\Delta t \sum_{K\in \Omega_R}\sum_{L\in N(K)}
    \frac{m(\sigma_{K,L})}{d(K,L)}\\
&\qquad \qquad \qquad \qquad \qquad 
\times (u^{n+1}_{1,L}-u^{n+1}_{1,K})(\varphi_1(t^n,x_{L})-\varphi_1(t^n,x_K)).
  \end{split}
\end{equation*}
Hence
\begin{equation*}
  \begin{split}
    T_2+T^*_2=&
\frac{1}{2} \sn a_1\Bigl(\sum_{K_0 \in \Om_h}u^{n}_{1,K_0}\Bigl)
\sum_{K\in \Omega_R} \sum_{L\in N(K)}
m(\sigma_{K,L}) (u^{n+1}_{1,L}-u^{n+1}_{1,K})\\
&\qquad\times \Biggl(\int_{t_n}^{t_{n+1}}
\frac{\varphi_1(t^n,x_{L})-\varphi_1(t^n,x_K)}{d(K,L)}\\
&\qquad \qquad \qquad
-\frac{1}{m(\sigma_{K,L})}
\int_{t_n}^{t_{n+1}}\int_{\sigma_{K,L}}\Grad \varphi_1\cdot\eta_{K,L}d\gamma\Biggl).
  \end{split}
\end{equation*}
Since the straight line $(x_K,x_{L})$ is orthogonal to $\sigma_{K,L}$, we have
$$
x_K-x_{L}= d(K,L)\eta_{K,L}. 
$$
This implies from the regularity of $\varphi_1$ that 
$$
\frac{\varphi_1(t^n,x_{L})-\varphi_1(t^n,x_K)}{d(K,L)}\equiv 
\Grad \varphi_1(t^n,x)\cdot \eta_{K,L}\text{ with $x$ between $x_K$ and $x_L$},
$$
and so
\begin{equation}\label{est:fv:dif:1}
  \begin{split}
&\Biggl|\int_{t_n}^{t_{n+1}}
\frac{\varphi_1(t^n,x_{L})-\varphi_1(t^n,x_K)}{d(K,L)}\\
&\qquad \qquad -\frac{1}{m(\sigma_{K,L})}
\int_{t_n}^{t_{n+1}}\int_{\sigma_{K,L}}\Grad \varphi_1\cdot\eta_{K,L}d\gamma \Biggl|
\le C \Delta t \, h,
  \end{split}
\end{equation}
for some constant $C>0$. Using \eqref{est:fv:dif:1} and \eqref{est:H1-norm}, we deduce
$$
\lim_{h\to 0} T_2=-\int_0^T a_1\Bigl(\int_\Om u_1 \dx\Bigl)
\int_\Om u_1 \Delta \varphi_1\dx \dt
=\int_0^T a_1\Bigl(\int_\Om u_1 \dx\Bigl) \int_\Om \Grad u_1 \cdot \Grad \varphi_1\dx \dt.
$$

Now, we show that
$$
\lim_{h \to 0}T_3=\int_0^T\int_\Om
\Bigl(\sigma(u_{1},u_{2},u_{3})
+\mu u_{1}\Bigl)\varphi_1 \dx \dt.
$$
For this purpose, we introduce
$$
T_{3,1}:=\sn \sum_{K\in \Omega_R} \Bigl(\sigma(u^{n+1}_{1,K},u^{n+1}_{2,K},u^{n+1}_{3,K})
+\mu u^{n+1}_{1,K}\Bigl)\int_{t^n}^{t^{n+1}}\int_K
\Bigl(\varphi_1(t^n,x_K)-\varphi_1(t,x)\Bigl)\dx \dt
$$
and
\begin{equation*}
  \begin{split}
T_{3,2}:=&\sn \sum_{K\in \Omega_R} \int_{t^n}^{t^{n+1}}\int_K
\Bigl(\sigma(u^{n+1}_{1,K},u^{n+1}_{2,K},u^{n+1}_{3,K})-\sigma(u_{1},u_{2},u_{3})
\Bigl)\varphi_1(t,x)\dx \dt\\
&\qquad \qquad +\sn \sum_{K\in \Omega_R} \int_{t^n}^{t^{n+1}}\int_K
\mu(u^{n+1}_{1,K}-u_{1})\varphi_1(t,x)\dx \dt.
  \end{split}
\end{equation*}
We have for all $x \in K$ and $t \in [t^n,t^{n+1}]$ that
\begin{equation}\label{equat-conv:1}
\abs{\varphi_1(t^n,x_K)-\varphi_1(t,x)}\le C (\Delta t+h),
\end{equation}
and thus, thanks to \eqref{est:L2-norm} and \eqref{est:H1-norm},
\begin{equation*}
\abs{T_{3,1}}\le C (\Delta t+h)\sn \Delta t \sum_{K\in \Omega_R} m(K)
\Bigl(\sigma(u^{n+1}_{1,K},u^{n+1}_{2,K},u^{n+1}_{3,K})+\mu u^{n+1}_{1,K}\Bigl)
\le C (\Delta t+h).
\end{equation*}
Hence, $T_{3,1}\to 0$ as $h\to 0$. We also have
\begin{equation*}
\abs{T_{3,2}}\le C \int_{0}^{T}\int_\Om
\Bigl(\abs{\sigma(u_{1,h},u_{2,h},u_{3,h})-\sigma(u_{1},u_{2},u_{3})}+\abs{u_{1,h}-u_1}\Bigl)\dx \dt.
\end{equation*}
Therefore from this and \eqref{limit-strong}, we deduce
$$
\text{$\abs{T_{3,2}}$ tends to zero as $h\to 0$}.
$$
This concludes the proof of \eqref{eq1:def1}. 
Reasoning along the same lines as above, we conclude that also 
\eqref{eq2:def1} and \eqref{eq3:def1} hold.

\section{Numerical Examples} \label{num-ex}


\subsection{Example 1. SIR model simulations.}
In this section we consider a sample square domain
$\Omega=(0,1)\times(0,1)$ and we show the behavior 
of the solution for different models of nonlocal functions
$a_i$, $i=1,2,3$.

\begin{figure}[h]
\begin{tabular}{cc}
(a) & (b) \\
\epsfig{file= 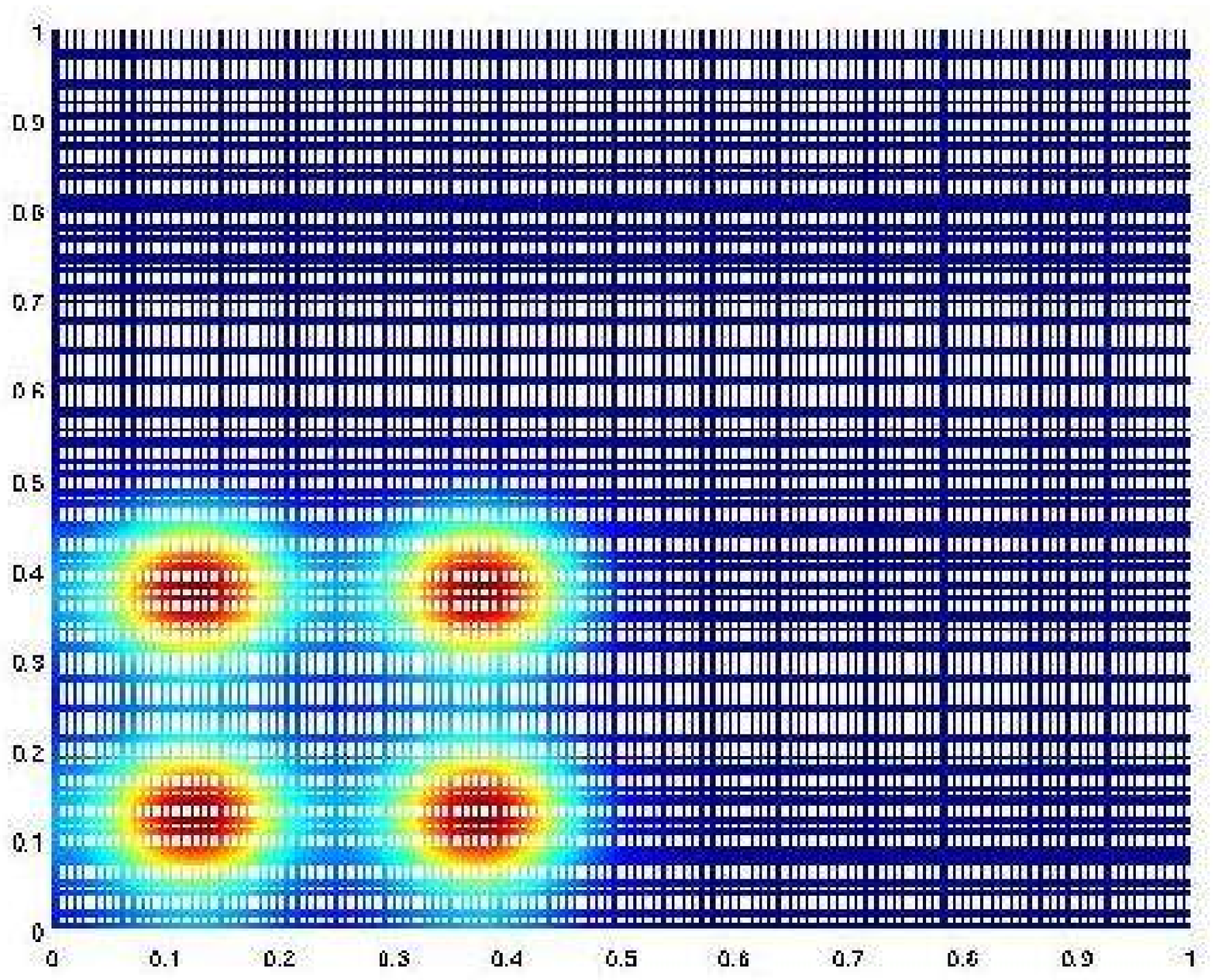,height=0.4\textwidth,width=0.4\textwidth} &
\epsfig{file= 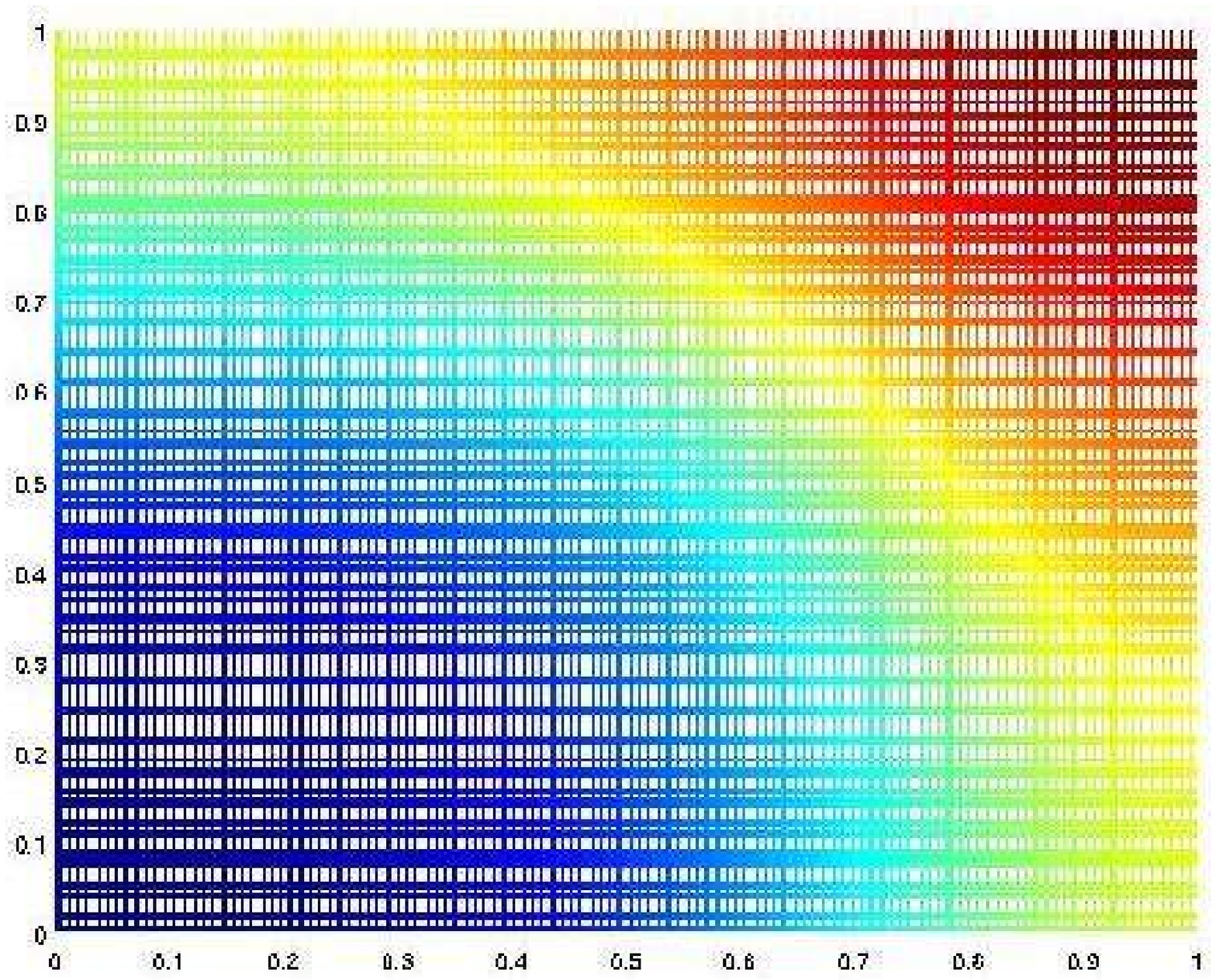,height=0.4\textwidth,width=0.4\textwidth} \\
(c) & (d) \\
\epsfig{file= 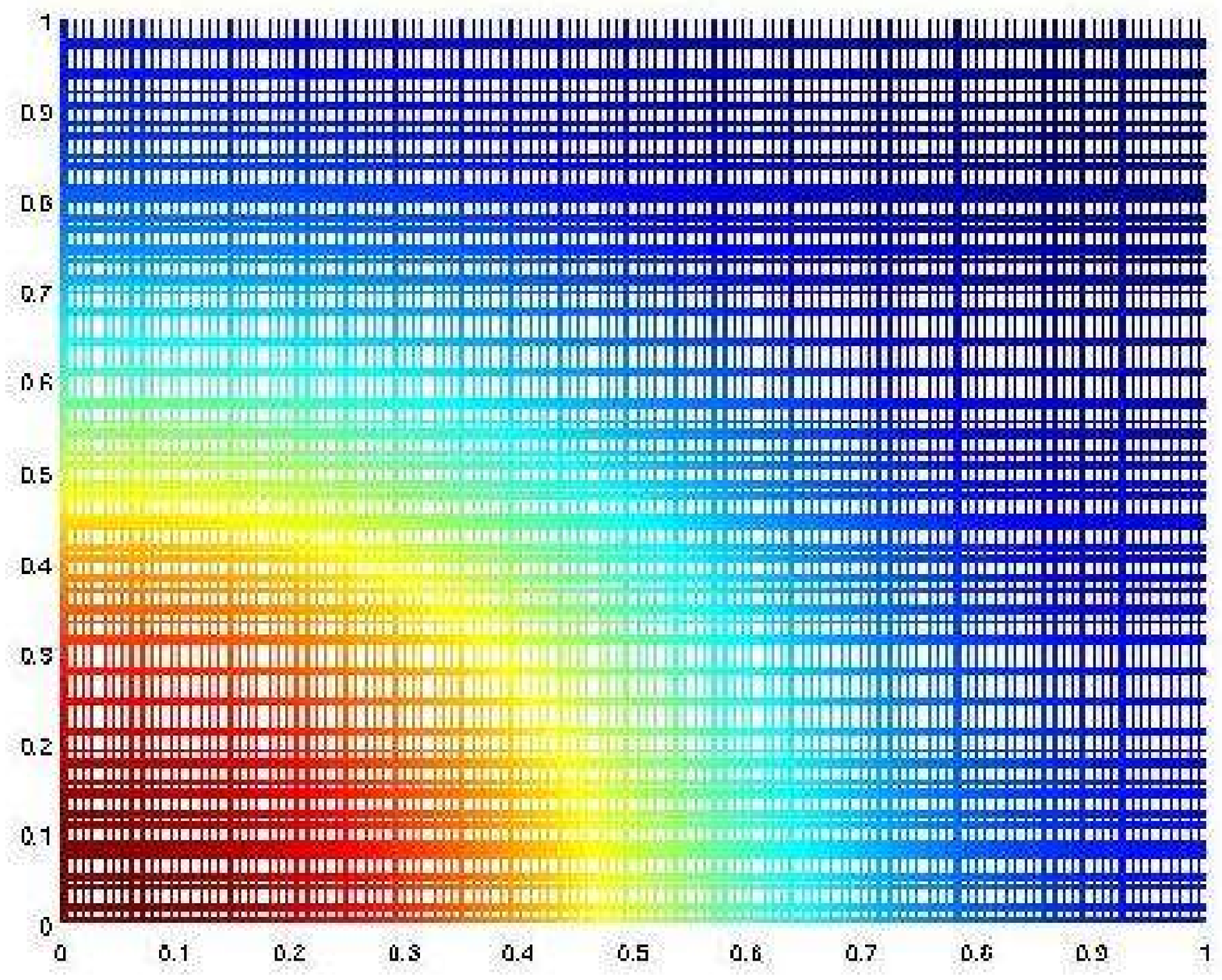,height=0.4\textwidth,width=0.4\textwidth} &
\epsfig{file= 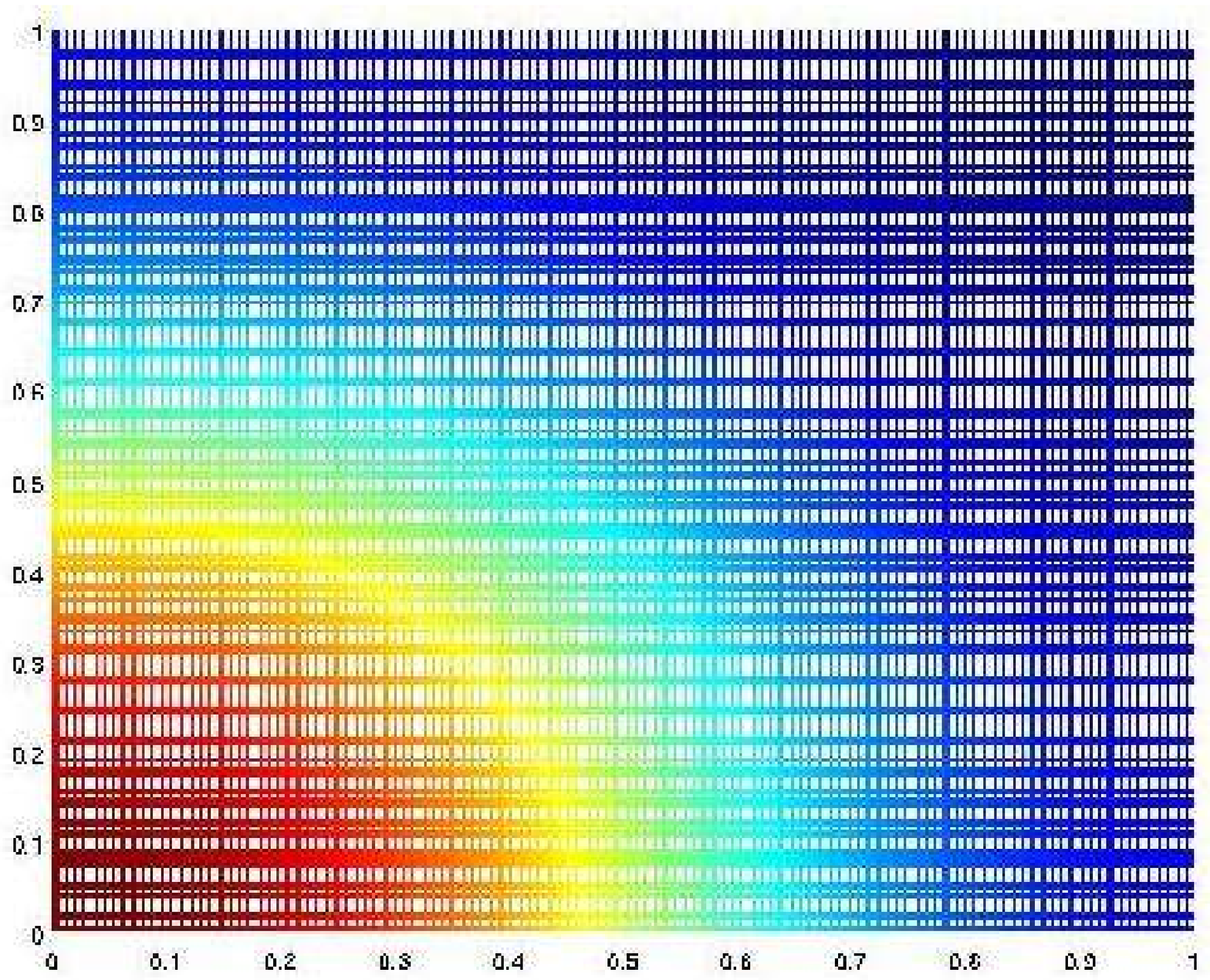,height=0.4\textwidth,width=0.4\textwidth}
\end{tabular}
\caption{SIR model with $a_1=a_2=a_3=1/10$:\label{fig:const_diff}
(a) Beginning infected population ($t=0.025\;sec.$); 
(b) Susceptible population ($t=0.5\;sec.$);
(c) Infected population ($t=0.5\;sec.$);
(d) Recovery population ($t=0.5\;sec.$).
}
\end{figure}


\begin{figure}[h]
\begin{tabular}{cc}
(a) & (b) \\
\epsfig{file= 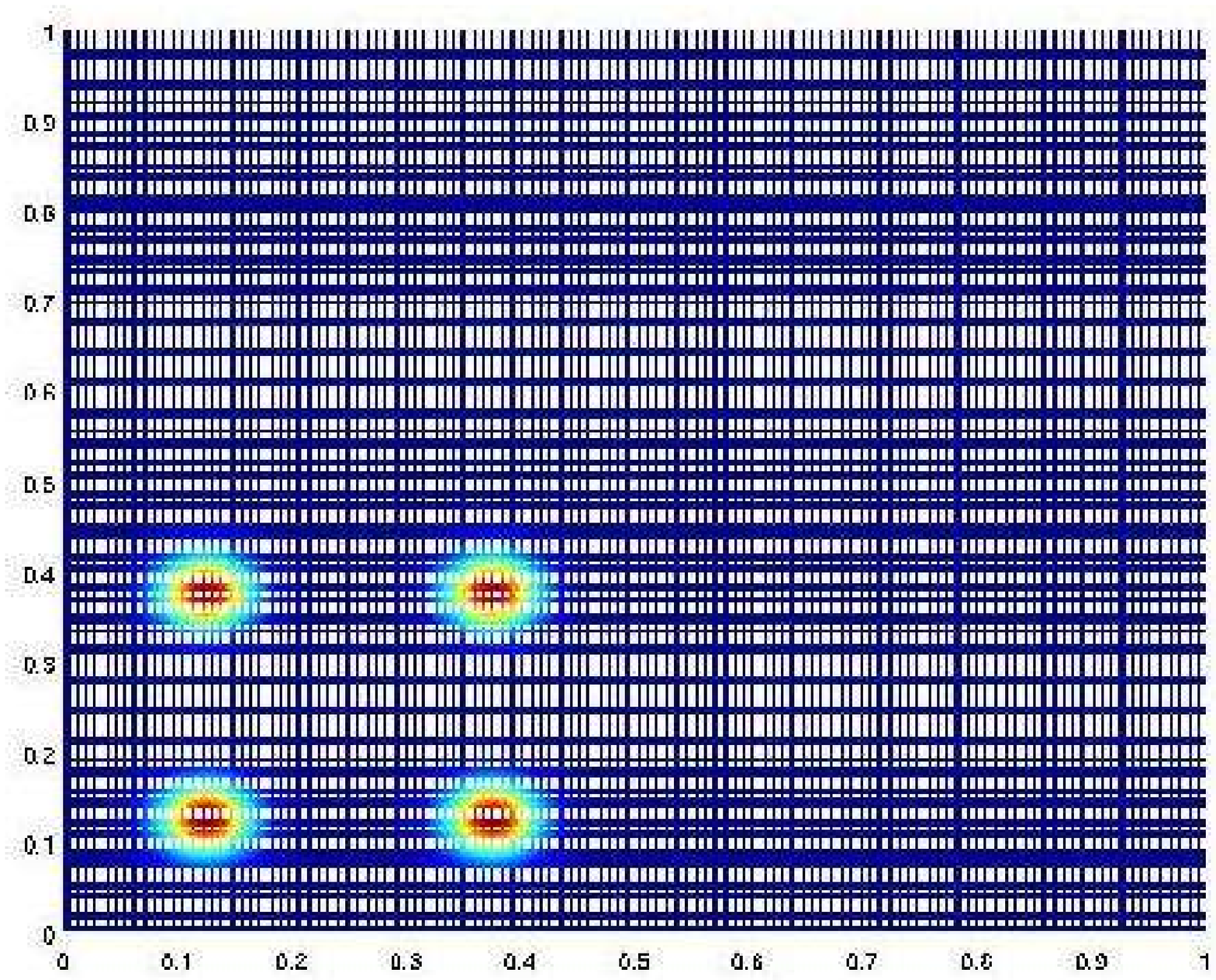,height=0.4\textwidth,width=0.4\textwidth} &
\epsfig{file= 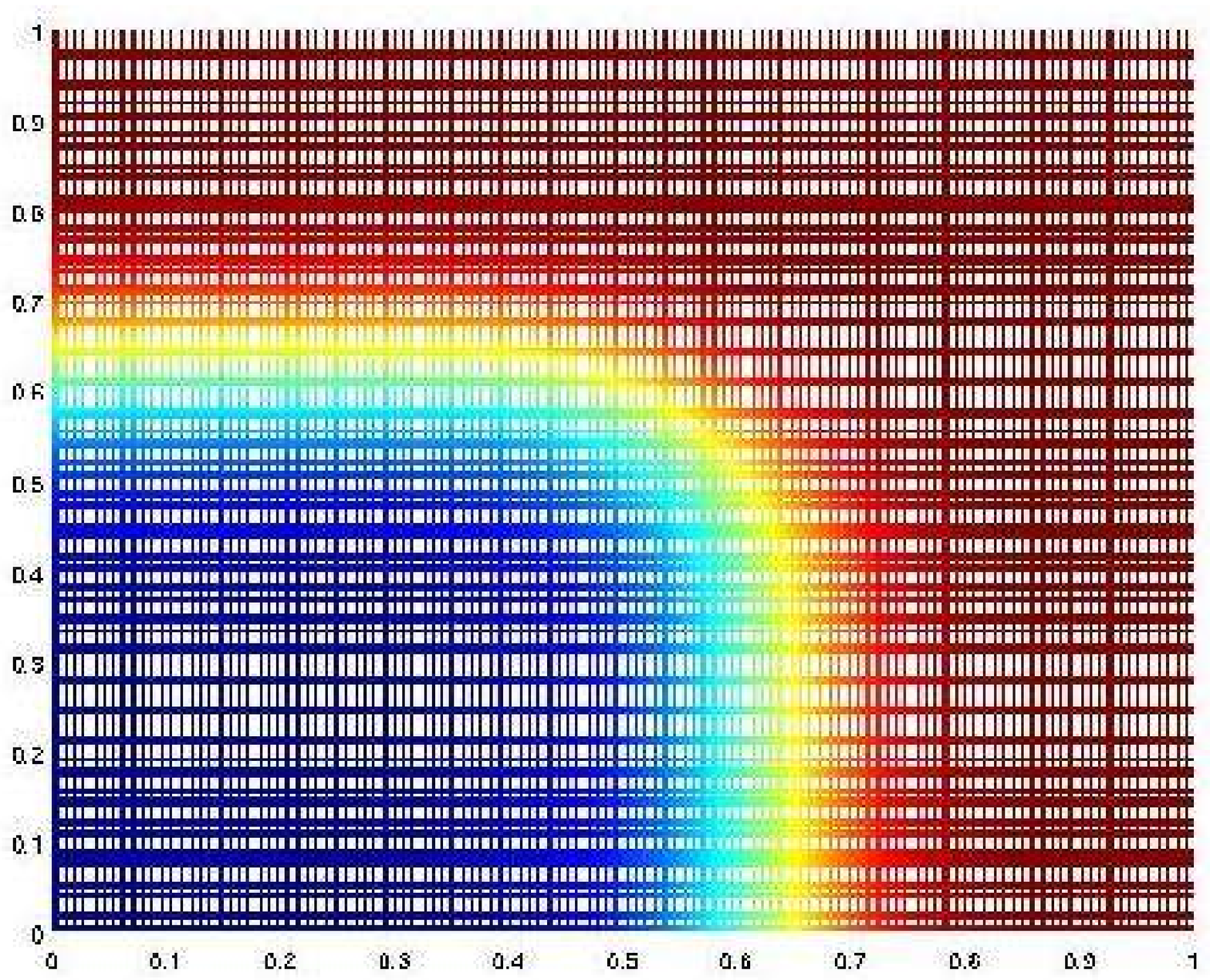,height=0.4\textwidth,width=0.4\textwidth} \\
(c) & (d) \\
\epsfig{file= 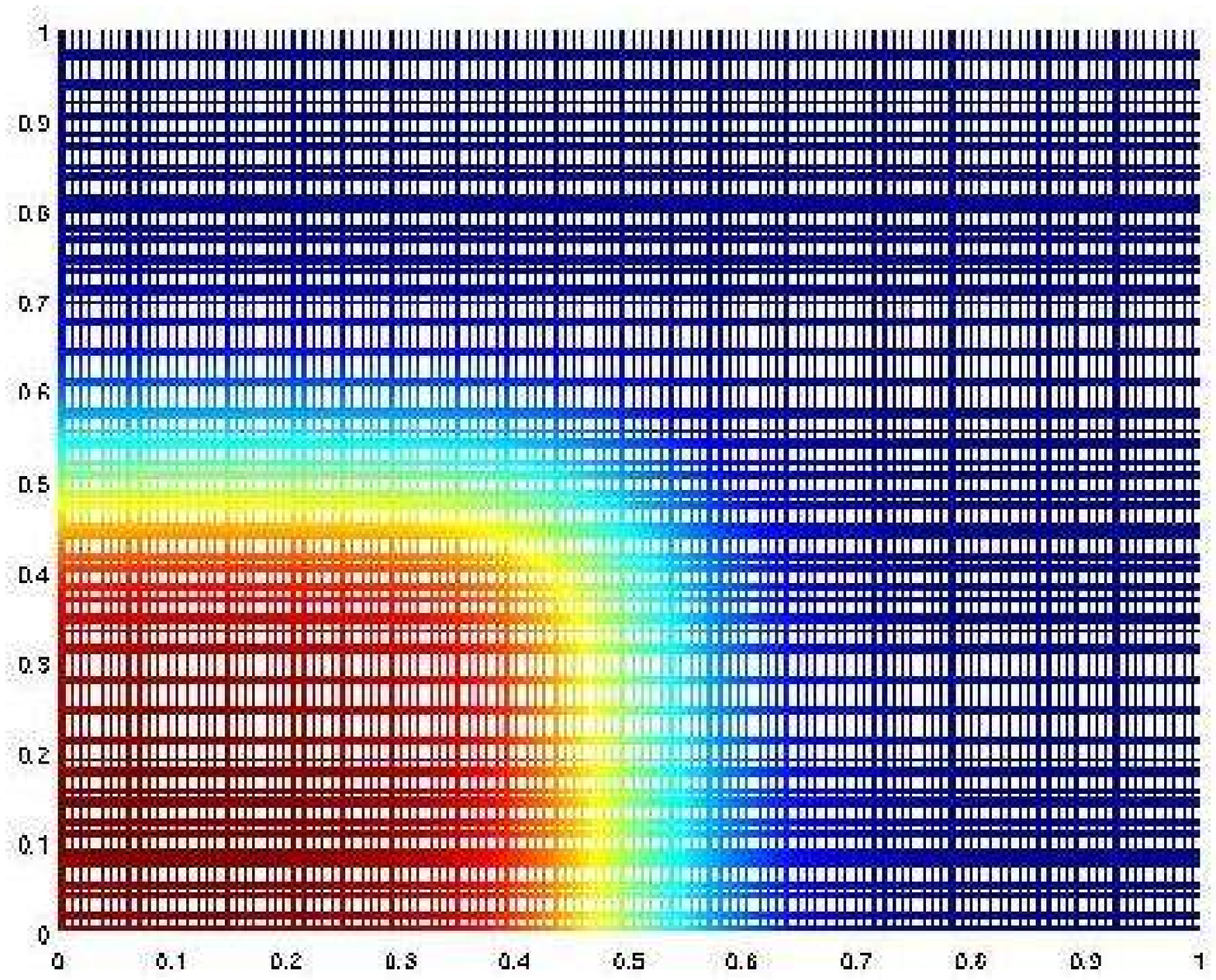,height=0.4\textwidth,width=0.4\textwidth} &
\epsfig{file= 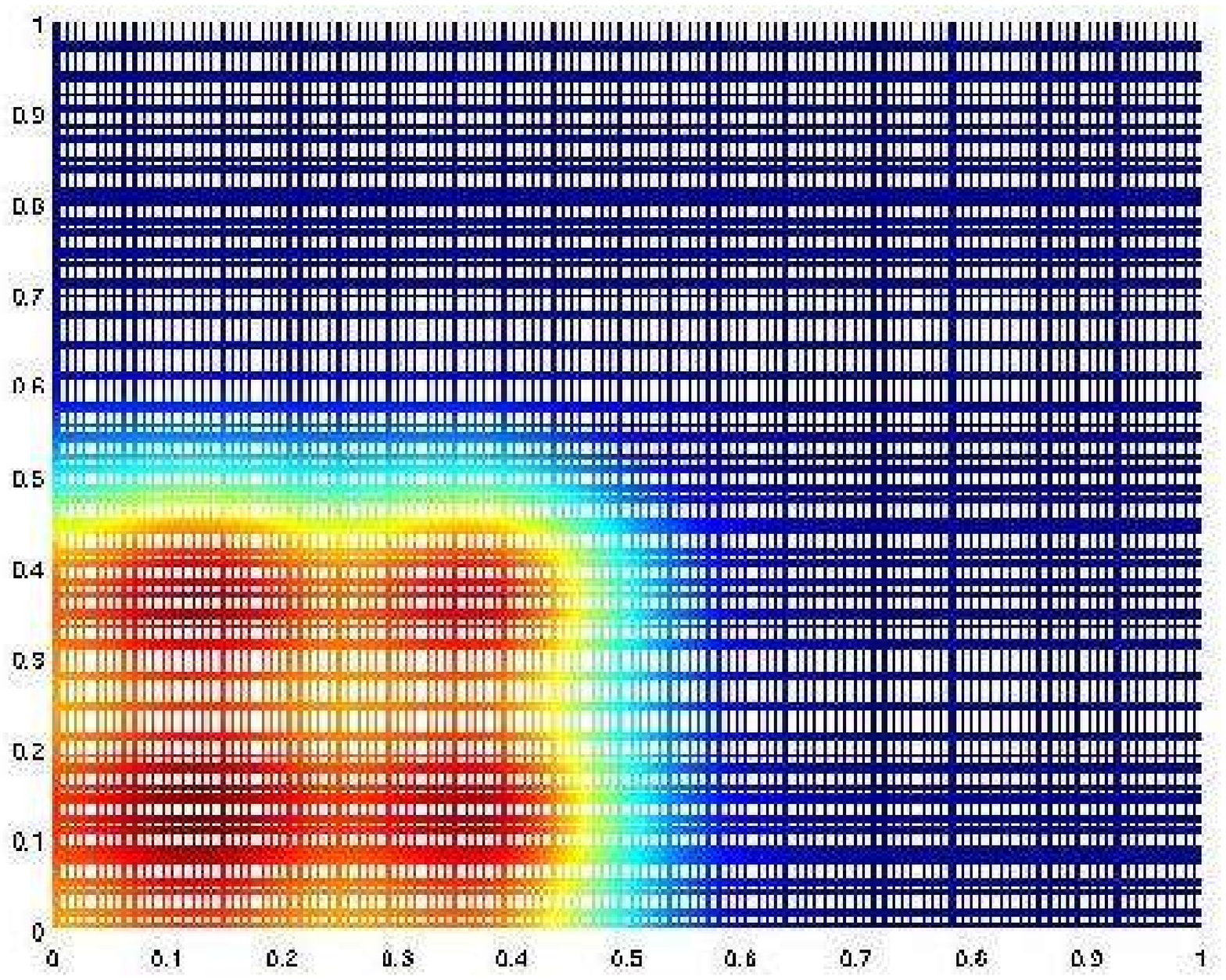,height=0.4\textwidth,width=0.4\textwidth}
\end{tabular}
\caption{SIR model with $a_1=a_2=a_3=s/10$:\label{fig:prop_diff}
(a) Beginning infected population; 
(b) Susceptible population;
(c) Infected population;
(d) Recovery population.
}
\end{figure}

We consider here a uniform mesh given by a Cartesian grid 
with $N_x \times N_y$ control volumes
and choosing $N_x=N_y=300$ for all simulations. 
Obviously, it is possible to consider also unstructured meshes,
but we will confine here to uniform mesh for simplicity of the 
simulated models.
The discretization in time is given by
 $N_t=100$ time steps for $T=0.5$. That is,
$\delta t=T/N_t$ and  $m(K)=1/(N_x N_y)$.
The parameter of the SIR model are given by
$\alpha=2.0$, $\mu=0.01$ and $\gamma=1.0$.


The initial condition are given by
\begin{eqnarray*}
\left\{
\begin{array}{l}
u_{1,0}(x,y)=\varepsilon_0;\\ 
u_{2,0}(x,y)=B\displaystyle \sum_{j=1}^5{\rm sech}(\beta(x-x_j)){\rm sech}(\beta(y-y_j));\\ 
u_{3,0}(x,y)=0.0.
\end{array}
\right.
\end{eqnarray*}
with
$\varepsilon_0=0.01$, $B=5000$, $\beta=2000$, $(x_1,y_1)=(0.25.0.25)$,
$(x_2,y_2)=(0.125.0.125)$, $(x_3,y_3)=(0.125.0.375)$,
$(x_4,y_4)=(0.375.0.125)$ and $(x_5,y_5)=(0.375.0.375)$.
These initial conditions represent one hand to a susceptible population initially with a low density and constant throughout the domain. 
On the other hand, 5 pockets of high  density infected population are located in the quadrant $[0,1/2]\times[0,1/2]$
which will diffuse the epidemic desease on the rest of the domain. Finally, we assume that there is no initially presence of recovery population.

\begin{figure}[h]
\epsfig{file= 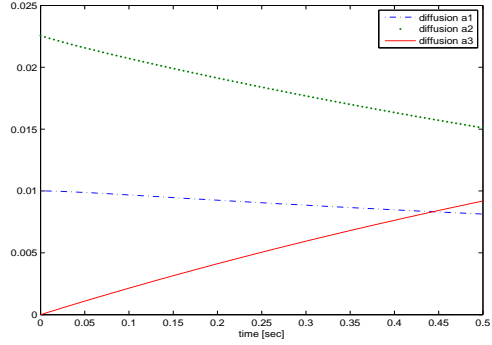,height=0.4\textwidth,width=0.6\textwidth}
\caption{Example 1: evolution in time of the diffusion terms 
$\displaystyle a_i\left(\int_\Omega u_i\,dx\right)$,
with $i=1,2,3$.\label{fig:diff_ex1}
}
\end{figure}


In a first simulation we consider a model which the diffusion rates do not depend on  the population,
that is $a_i$, $i=1,2,3$ are constant and equal to $0.1$ (see Figure \ref{fig:const_diff}).
On the other hand, in figure \ref{fig:prop_diff} we observe another simulations with 
a nonlocal diffusion where
the three diffusion rates are linearly proportional to
the total mass of each population ($a_1(s)=a_2(s)=a_3(s)=0.1\,s$). 
Outside the nonlocal diffusion, we consider the same parameters for both simulations.
We remark that second simulation with nonlocal diffusion violate the assumption \eqref{ass-ai}. In fact, they correspond 
to degenerated parabolic cases. In order to ensure the convergence of our numerical example, we 
replace in a practical way, the diffusion rate coefficients by
\begin{equation}
\label{trunc}
\widetilde{a}_i(s) = 
\left\{
\begin{array}{ll}
M& \hbox{ if } s>M\\
s& \hbox{ if } \varepsilon\leqslant s \leqslant M\\
\varepsilon& \hbox{ if } s<\varepsilon,
\end{array}
\right.
\end{equation}
with $M=10^4$ and $\varepsilon=10^{-4}$.
Figures \ref{fig:const_diff} and \ref{fig:prop_diff} represent the simulation at time
$t=0.025$ for the localized infected population (see picture (a) of both figures)
and the simulation at time $t=0.5$ for the different populations (see pictures (b), (c) and (d)). 
We remark that different behaviours occur
for the susceptible population, observing a first extinsion of this population in a different zones of the square $[0,1]\times[0,1]$ for each case.
Finally, we observe in Figure \ref{fig:diff_ex1}, the
 evolution in time of the diffusion
$\displaystyle a_i\left(\int_\Omega u_i\,dx\right)$, with $i=1,2,3$ 
for the simulation with the nonlocal diffusion.

Unfortunately, it is not possible to observe patterns formation for this SIR-diffusion model like some 
results obtained  in \cite{liu}. That is
because  the infected population have his equilibrium in $u_2=0.0$. For this reason, it is not difficult 
to compute the Jacobian matrix for the SIR reaction term and to prove that 
the diffusion coefficients do not affect the sign of the real part of the eigenvalues.
In other, words, the Turing space is always empty in this case (see \cite{turing}), and the populations tend to their 
constant equilibrium state for long times.

\subsection{Example 2. An epidemic model of SARS with Patterns Formation.}

Now, we consider the following modification of the SIR model 
\begin{equation}
	\label{S1_liu}
	\begin{cases}
		&\displaystyle \pt u_1 
		-a_1\Bigl(\int_\Om u_1\dx\Bigl) \Delta u_1= A -\sigma(u_1,u_2,u_3)-\mu u_1,\\ 
		&\displaystyle \pt u_2 
		-a_2\Bigl(\int_\Om u_2\dx\Bigl) \Delta u_2 = \sigma(u_1,u_2,u_3)-\gamma u_2-\mu u_2 - \mathcal{H}(u_2),\\ 
		&\displaystyle \pt u_3 
		-a_3\Bigl(\int_\Om u_3\dx\Bigl)\Delta u_3= \gamma u_2 +\mathcal{H}(u_2)-\mu u_3.
	\end{cases}	
\end{equation}
In this modified SIR model,
$A$ is the recruitment rate of the population (such as growth rate of average population size, a recover becomes an susceptible, immigrant and so on), and $\mathcal{H}(u_2)$ is the removal rate of infective individuals due to the treatment.
We suppose that the treated infectives becomes recovered when they are treated in treatment, and
$$
\mathcal{H}(v)=
\left\{
\begin{array}{lll}
r&\mbox{ if } &v>0,\\
0&\mbox{ otherwise },
\end{array}
\right.
$$
where $r>0$ is constant and represents the capacity of treatment for infectives. 
We consider the same incidence term $\sigma$ defined in the introduction \eqref{def-sigma-1},
and the natural death rate $\mu$ is included for the recovery population. 
The detail about this epidemic model
can be found in \cite{liu,wang}
where the authors adopt a bilinear incidence rate more simple than \eqref{def-sigma-1}.
One of the application to consider this model is that it
 supposes that the capacity for the treatment of a disease in a community
is a constant $r$, in order to  use the maximal treatment capacity to cure
or isolate infectives so that the disease is eradicated \cite{wang}. 
It can be used for example for
 mathematical model to simulate the SARS outbreak in Beijing \cite{wang2}. 

The approximation  of this modified SIR model using Finite Volume
is very similar to \eqref{S1-discr}-\eqref{S3-discr}.
More precisely we have
\begin{equation}\label{S1-discr_liu}
     \begin{split}
     &m(K)\frac{u^{n+1}_{1,K}-u^{n}_{1,K}}{\Delta t}
	 -a_1\Bigl(\sum_{K_0 \in \Om_h}u^{n}_{1,K_0}\Bigl)\sum_{L \in N(K) }
     \frac{m(\sigma_{K,L})}{d(K,L)}(u^{n+1}_{1,L}-u^{n+1}_{1,K})
     \\ &\qquad \qquad \qquad 
     +m(K)\left(-A+\sigma(u^{n+1,+}_{1,K},u^{n+1,+}_{2,K},u^{n+1,+}_{3,K})+\mu u^{n+1}_{1,K}\right)=0,
     \end{split}
\end{equation}
\begin{equation}\label{S2-discr_liu}
     \begin{split}
     &m(K)\frac{u^{n+1}_{2,K}-u^{n}_{2,K}}{\Delta t}
	 -a_2\Bigl(\sum_{K_0 \in \Om_h}u^{n}_{2,K_0}\Bigl)\sum_{L \in N(K) }
     \frac{m(\sigma_{K,L})}{d(K,L)}(u^{n+1}_{2,L}-u^{n+1}_{2,K})
     \\ &\qquad  
     -m(K)\left(\sigma(u^{n,+}_{1,K},u^{n+1,+}_{2,K},u^{n+1,+}_{3,K})-(\gamma+\mu) u^{n+1}_{2,K} -\mathcal{H}(u^{n+1}_{2,K})\right)=0,
     \end{split}
\end{equation}
\begin{equation}\label{S3-discr_liu}
     \begin{split}
     &m(K)\frac{u^{n+1}_{3,K}-u^{n}_{3,K}}{\Delta t}
	 -a_3\Bigl(\sum_{K_0 \in \Om_h}u^{n}_{3,K_0}\Bigl)\sum_{L \in N(K) }
     \frac{m(\sigma_{K,L})}{d(K,L)}(u^{n+1}_{3,L}-u^{n+1}_{3,K})
     \\ &\qquad \qquad \qquad\qquad \qquad \qquad 
     -m(K)\left(\gamma u^{n}_{2,K}+\mathcal{H}(u^{n+1}_{2,K})-\mu u^{n+1}_{3,K}\right)=0.
     \end{split}
\end{equation}
The results of existence of solution of the Finite Volume scheme (Proposition \ref{prop:exist-fv}),
nonnegativity of the scheme (Lemma \ref{lem:nonnegativity}) and convergence main result
(Theorem \ref{theo1}) can be easy generalized using straightforward calculations to this modified SIR model 
\eqref{S1_liu} and his 
Finite Volume scheme \eqref{S1-discr_liu}-\eqref{S3-discr_liu}.

\subsubsection{Analysis of spatial patterns}

The equilibrium state of the system \eqref{S1_liu} is not exactly the same described in \cite{liu}
because we consider here a different incidence term than the given by Liu which takes the reduced form $\sigma=\alpha uv$. 
In our case, the two positive equilibrium points are given by $E_1=(u_1,v_1,w_1)$ and $E_2=(u_2,v_2,w_2)$ where
\begin{eqnarray}
\label{eq_v} v_{1,2}&=&\frac{A-r}{2\,R_0} -\frac{A}{2\,\alpha} \pm 
\frac{\sqrt{\left(r\,\alpha-A\,R_0-A\,\alpha)^2-4\,A^2\,R_0\,\alpha\right)}}{2\,\alpha\, R_0}\\
\label{eq_u} u_{1,2}&=& \frac{A-r-R_0\,v_{1,2}}{\mu},\qquad  w_{1,2}\ =\ {\frac {\gamma\,v_{1,2}+r}{\mu}},
\end{eqnarray}
with $R_0=\mu+\gamma$. The linear stability of the system \eqref{S1_liu} is obtained when the real part of
the eigenvalues of
the Jacobian matrix
\begin{equation*}
J=
\left[ \begin {array}{ccc} -{\frac {\alpha\,v}{u+v+w}}+{\frac {\alpha
\,uv}{ \left( u+v+w \right) ^{2}}}-\mu&-{\frac {\alpha\,u}{u+v+w}}+{
\frac {\alpha\,uv}{ \left( u+v+w \right) ^{2}}}&{\frac {\alpha\,uv}{
 \left( u+v+w \right) ^{2}}}\\\noalign{\medskip}{\frac {\alpha\,v}{u+v
+w}}-{\frac {\alpha\,uv}{ \left( u+v+w \right) ^{2}}}&{\frac {\alpha\,
u}{u+v+w}}-{\frac {\alpha\,uv}{ \left( u+v+w \right) ^{2}}}-\gamma-\mu
&-{\frac {\alpha\,uv}{ \left( u+v+w \right) ^{2}}}\\\noalign{\medskip}0
&\gamma&-\mu\end {array} \right] 
\end{equation*}
are negative. The eigenvalues of this Jacobian matrix are given by $\lambda_1=-\mu<0$ and $\lambda_2$ and $\lambda_3$ which
are root of
$$
\lambda^2+\left((2\mu+\gamma)+\alpha\frac{v-u}{u+v+w}\right)\lambda+\left((\mu+\gamma)\mu+\alpha\frac{(\mu+\gamma)v-\mu\,u}{u+v+w}\right)=0.
$$
Real part of $\lambda_2$ and $\lambda_3$ are negative if and only if all the coefficients of this polynomial function of degree 2 are positive,
that is
\begin{equation}
\label{stability_cond}
\frac{u+v+w}{\alpha}>\max\left\{\frac{u-v}{2\,\mu+\gamma}\;,\;\frac{u}{\mu+\gamma}-\frac{v}{\mu}\right\}.
\end{equation}
\begin{figure}[h]
\begin{tabular}{cc}
(ric) & (a) \\
\epsfig{file= 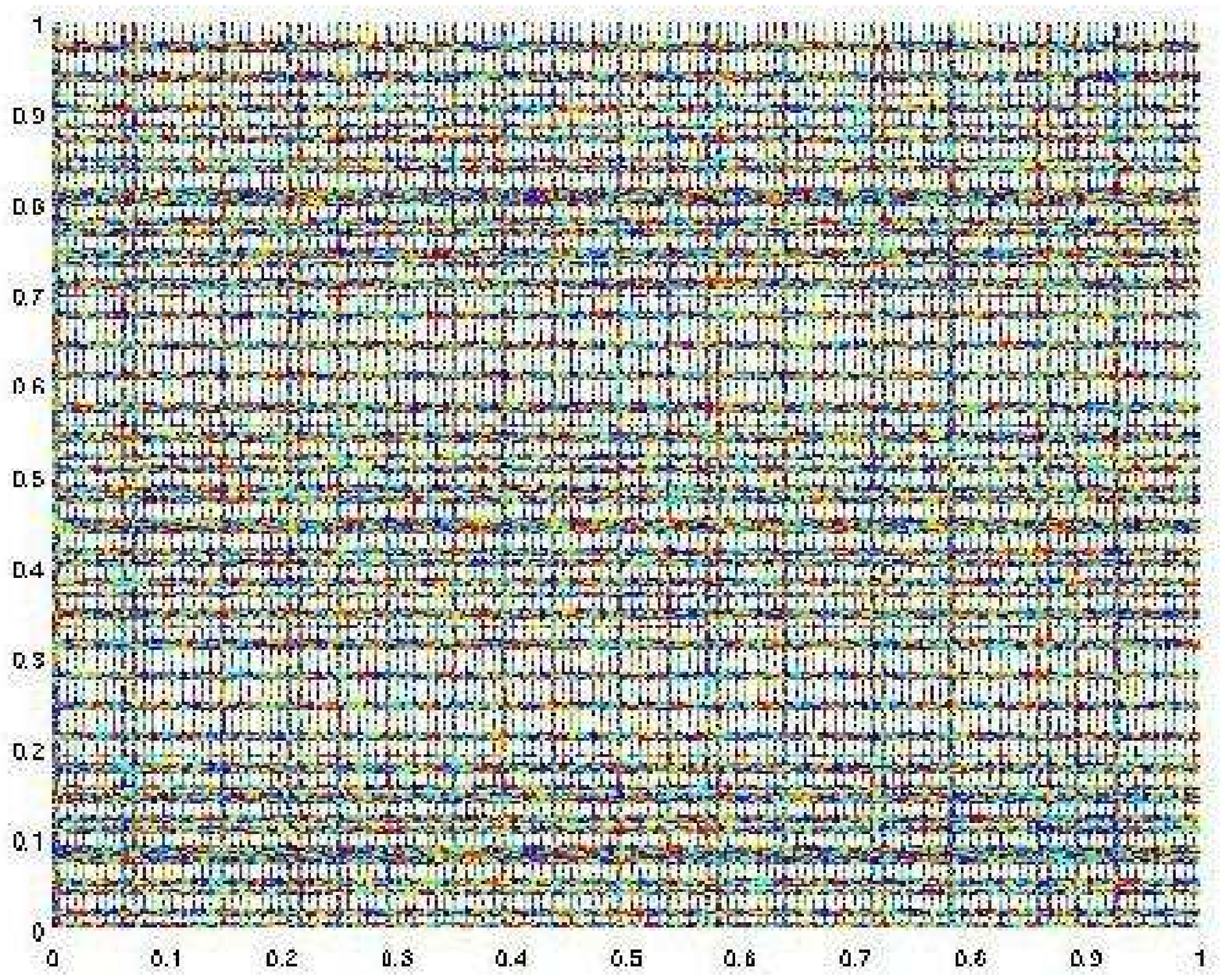,height=0.4\textwidth,width=0.4\textwidth} &
\epsfig{file= 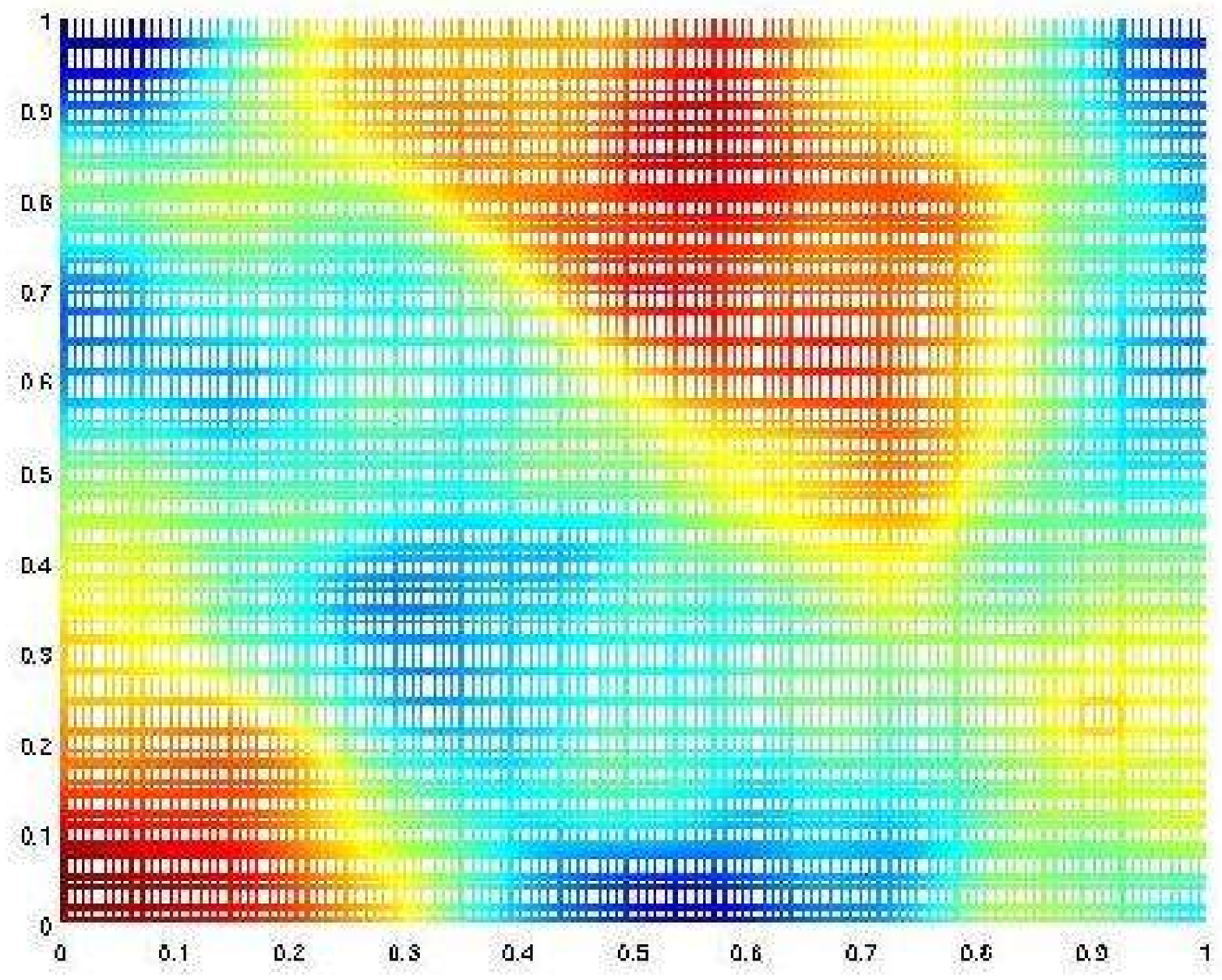,height=0.4\textwidth,width=0.4\textwidth} \\
(b) & (c) \\
\epsfig{file= 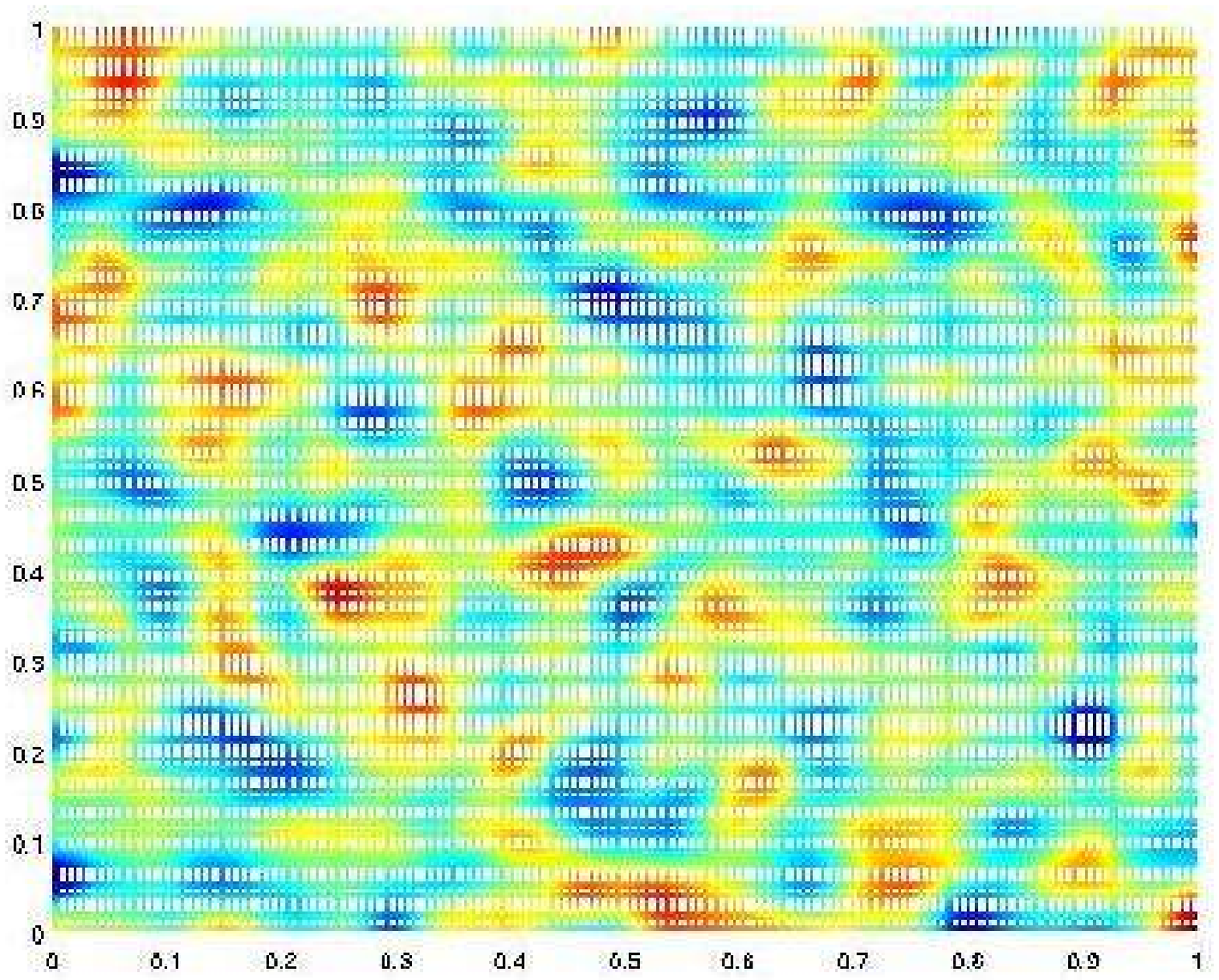,height=0.4\textwidth,width=0.4\textwidth} &
\epsfig{file= 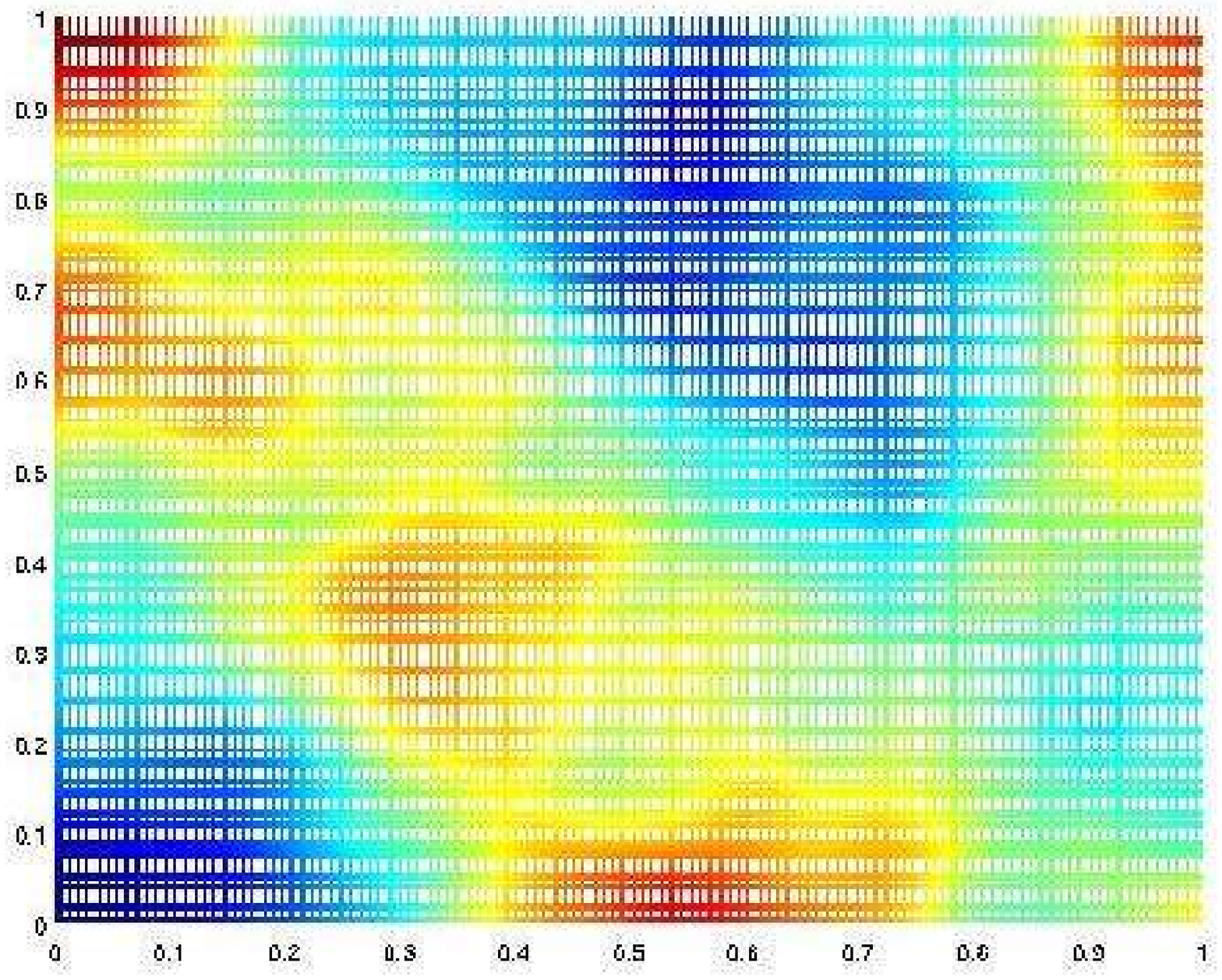,height=0.4\textwidth,width=0.4\textwidth}
\end{tabular}
\caption{Patter Formation for the SARS model (constant diffusion):\label{fig:SARS_const_diff}
(ric) Random initial condition; 
(a) Susceptible population;
(b) Infected population;
(c) Recovery population.
}
\end{figure}
\begin{figure}[h]
\begin{tabular}{cc}
\multicolumn{2}{c}{(a)}\\
\multicolumn{2}{c}{
\epsfig{file= 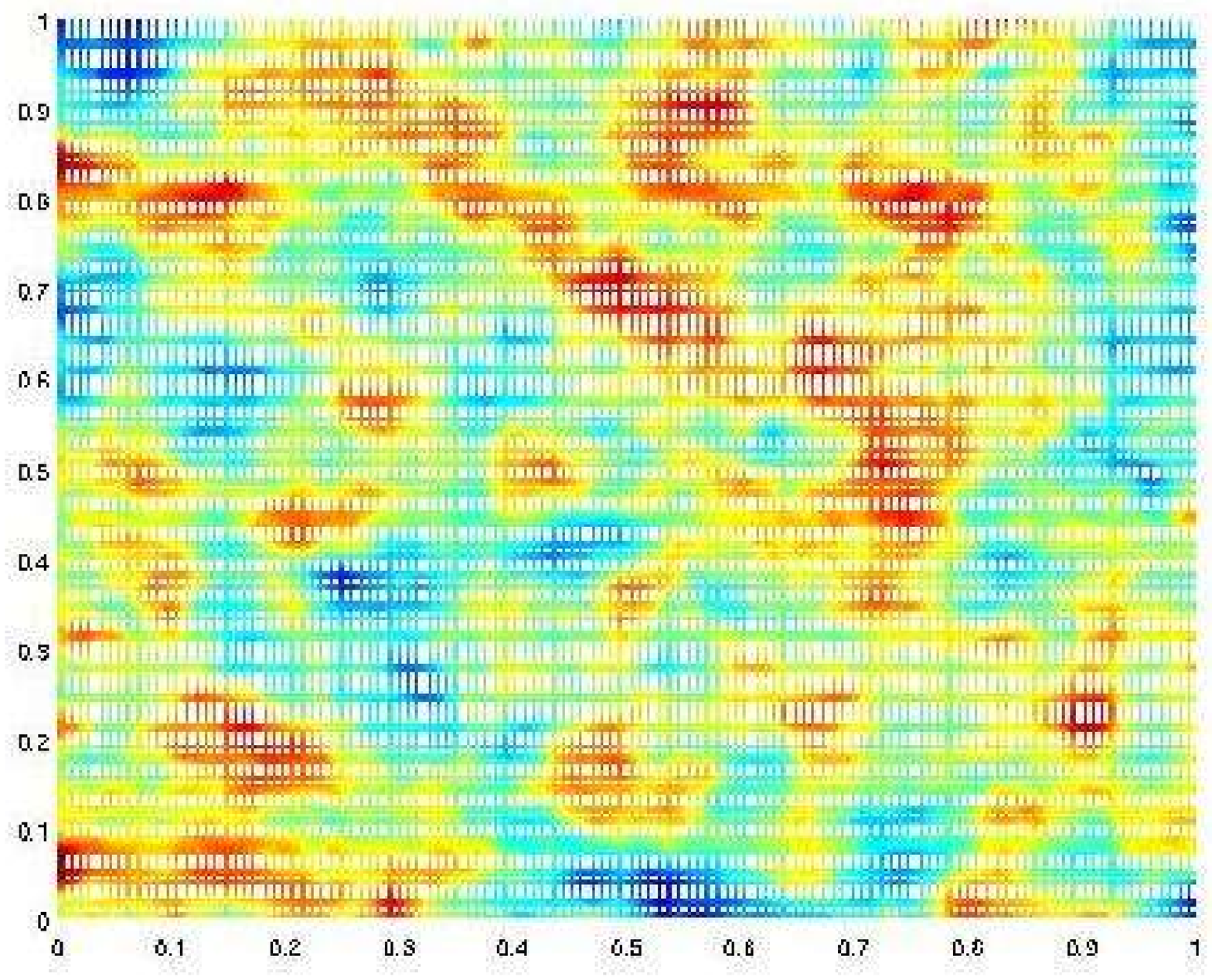,height=0.4\textwidth,width=0.4\textwidth} 
}
\\
(b) & (c) \\
\epsfig{file= 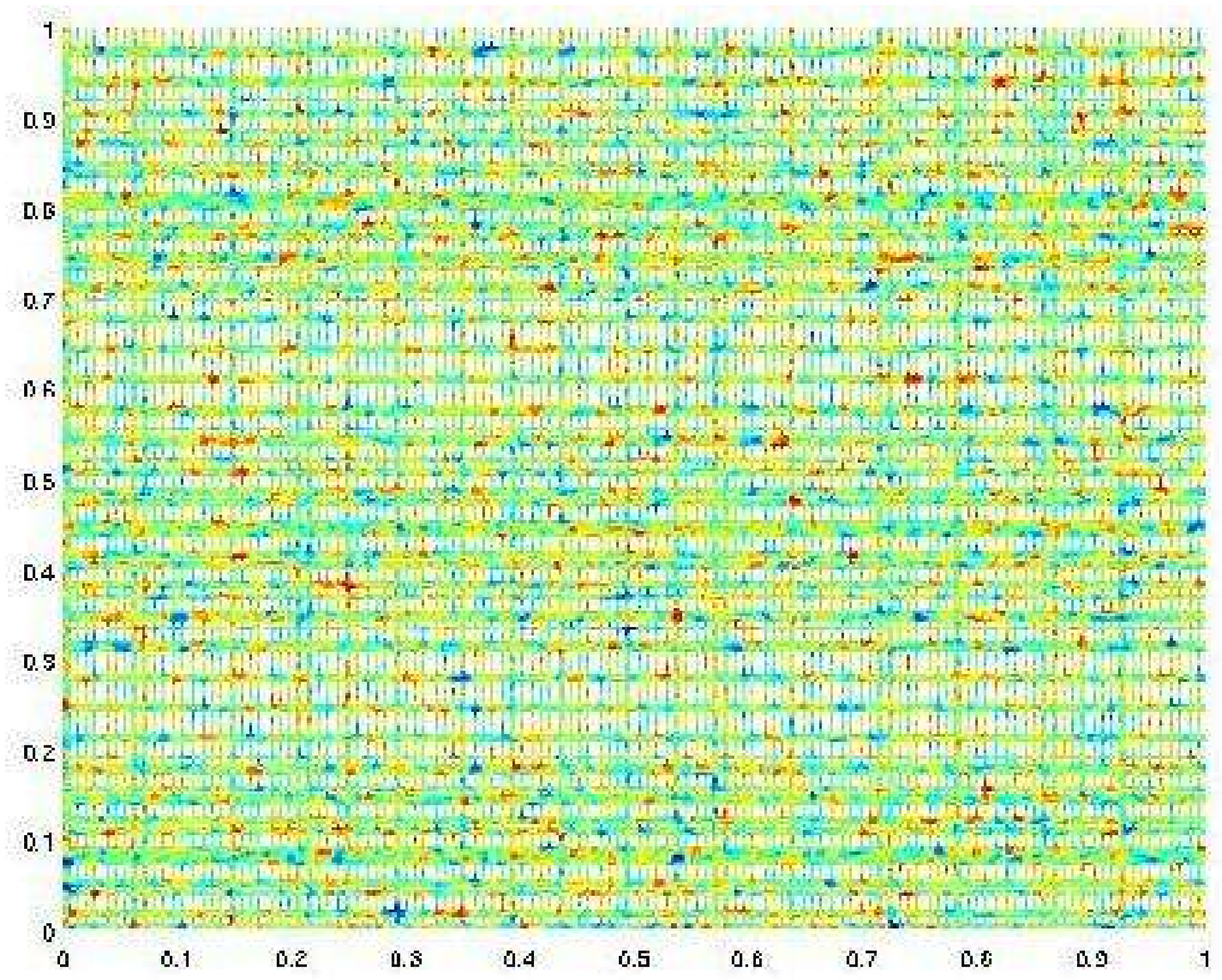,height=0.4\textwidth,width=0.4\textwidth} &
\epsfig{file= 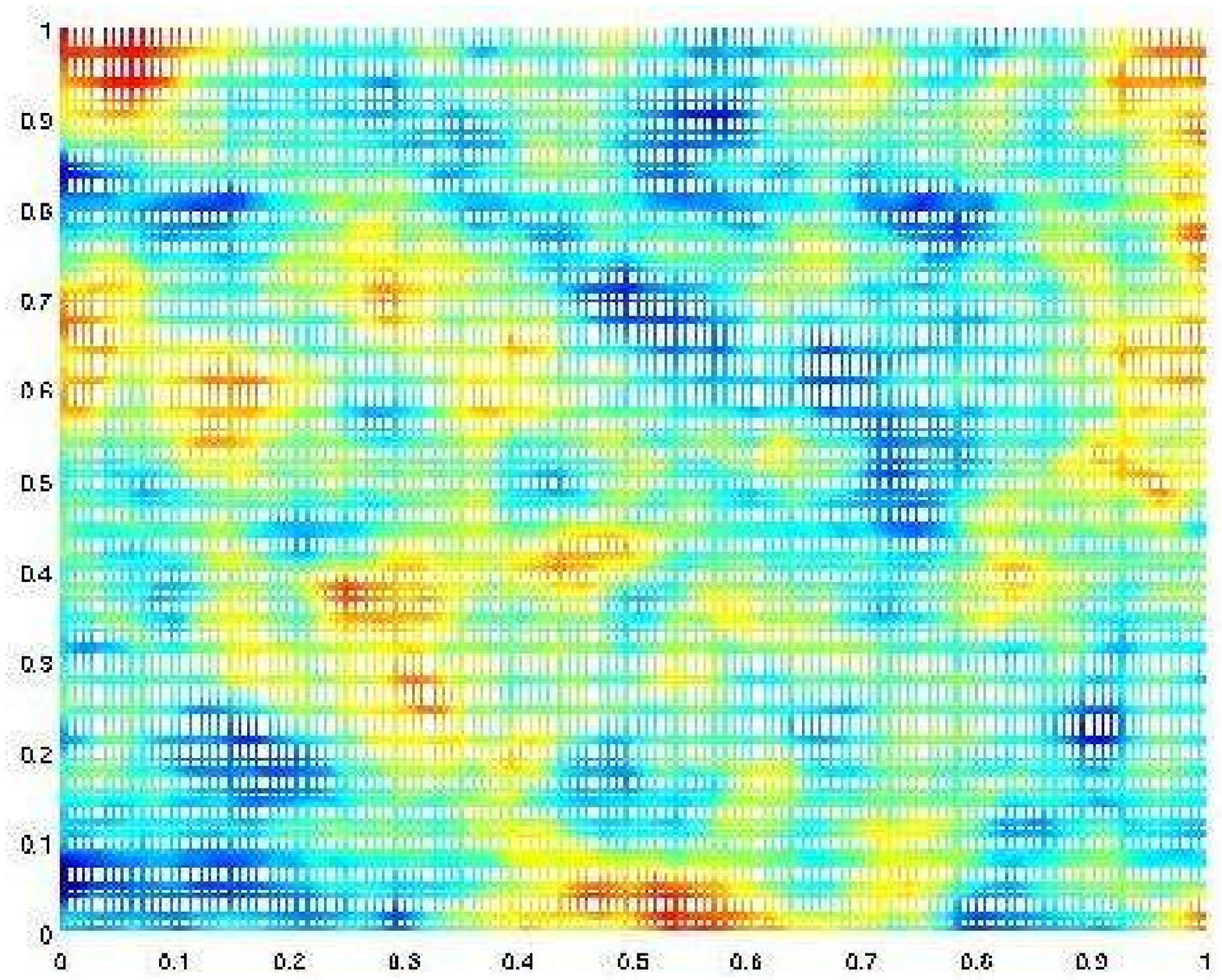,height=0.4\textwidth,width=0.4\textwidth}
\end{tabular}
\caption{Patter Formation for the SARS model (nonlocal diffusion):\label{fig:SARS_inv_prop_diff}
(a) Susceptible population;
(b) Infected population;
(c) Recovery population.
}
\end{figure}
\\
We consider in this example the parameters $A=3$, $\mu=0.3$, $\alpha=3.8$, $r=0.5$, and $\gamma=0.8$
(they are the same of the example of Liu \cite{liu}, except the parameter $\alpha$ of the incidence term).
Replacing this parameter in \eqref{eq_v}-\eqref{eq_u}, we have the equilibrium state 
\begin{eqnarray*}
E_1 &=& \left(7.217163781\;,\; 0.3044098832\;,\; 2.478426355 \right)\\
E_2 &=& \left(4.010906415\;,\; 1.178843705\;,\; 4.810249881 \right).
\end{eqnarray*}
It is easy to verify that the equilibrium point $E_2$ verify the linear stability condition \eqref{stability_cond}
and $E_1$ corresponds to a instable point.
Formation of spatial patterns 
results from the diffusion-induced instability when 
the real part of at least one of the eigenvalues of 
$J=\left(
\begin{array}{ccc}
d_1&0&0\\
0&d_2&0\\
0&0&d_3
\end{array}
\right)
$
is positive \cite{turing}. If we suppose that $d_1=d_3$, choosing the parameters above
and the the equilibrium point $E_2$, then Turing stabilities appear 
when
$(5.114590203+9\,d_1)\,d_2+3.289620038-2.20024467\,d1<0$.

\begin{figure}[h]
\begin{tabular}{cc}
\epsfig{file= 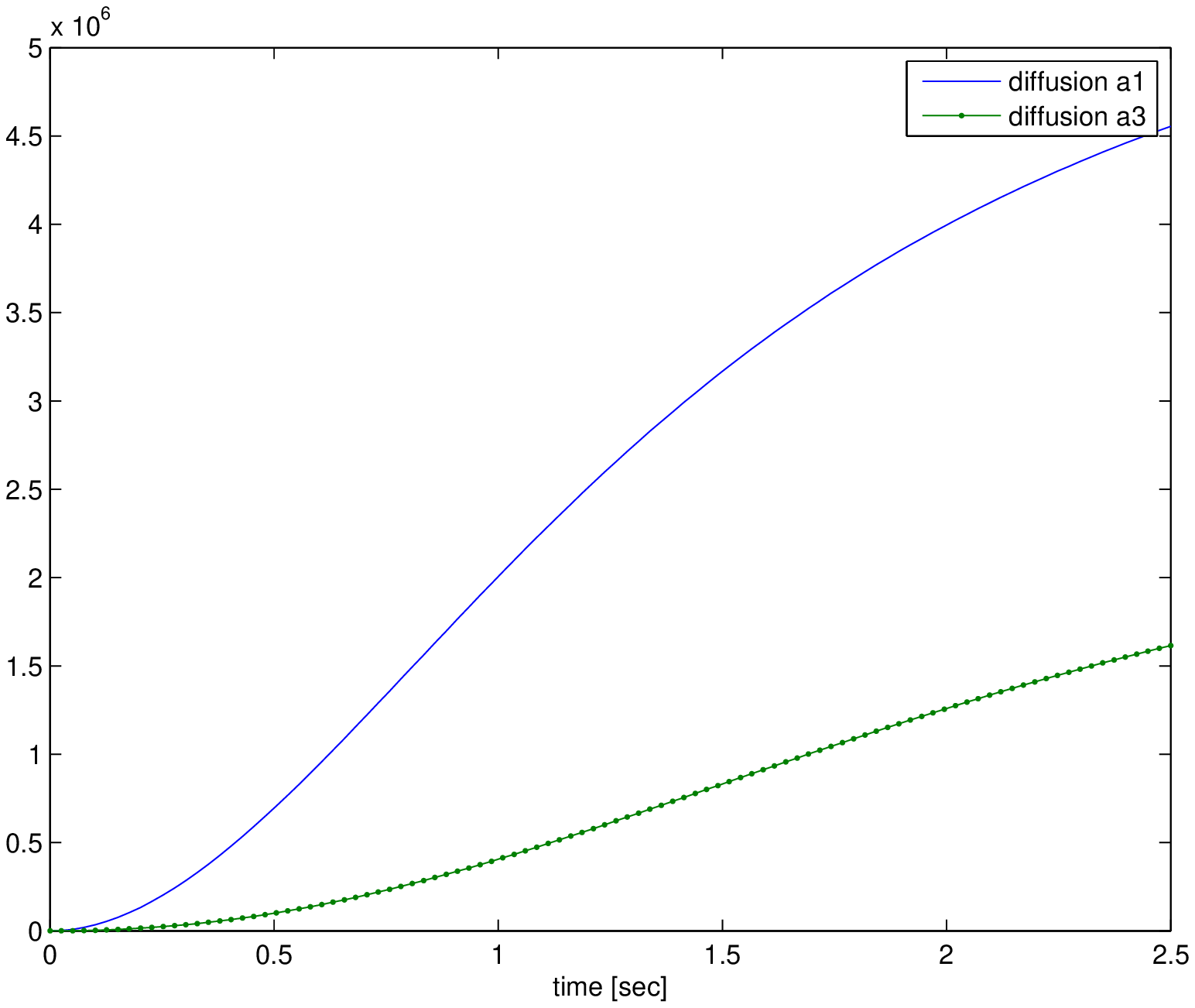,height=0.4\textwidth,width=0.5\textwidth}&
\epsfig{file= 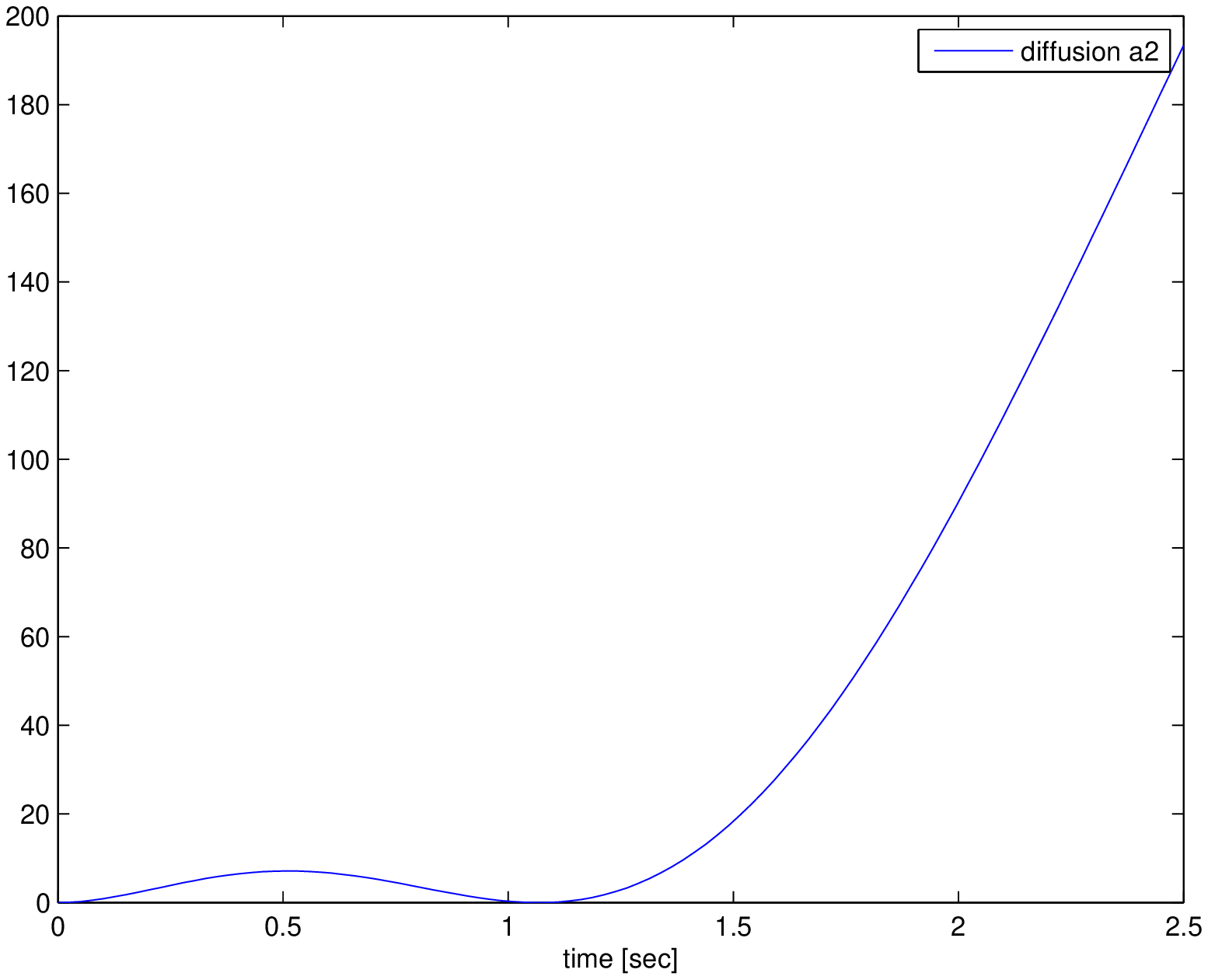,height=0.4\textwidth,width=0.5\textwidth}
\end{tabular}
\caption{Example 1: evolution in time of the diffusion terms 
$\displaystyle a_i\left(\int_\Omega u_i\,dx\right)$,
with $i=1,2,3$;
Left: Susceptible and recovery populations; Right:
Infected population.\label{fig:diff_ex2}
}
\end{figure}

\subsubsection{Simulations with local and nonlocal diffusion.}
 
Similar the above Example 1, we consider here two simulations, one with a constant diffusion and another one
with nonlocal diffusion. 
We take the same square domain of example 1, with 
 a Cartesian grid 
and choosing $N_x=N_y=300$ for both simulations. 
The discretization in time is given by
 $N_t=100$ time steps for $T=2.5$. That is,
$\delta t=T/N_t$ and  $m(K)=1/(N_x N_y)$.
The parameter of the SARS model are given by
$A=3$, $\mu=0.3$, $\alpha=3.8$, $r=0.5$, and $\gamma=0.8$.

In order to observe patterns formation, 
we consider the initial condition as follow:
\begin{eqnarray*}
\left\{
\begin{array}{l}
u_{1,0}(x,y)=\widetilde{u}_1 + \varepsilon_1 \omega_1;\\ 
u_{2,0}(x,y)=\widetilde{u}_2 + \varepsilon_2 \omega_2;\\ 
u_{3,0}(x,y)=\widetilde{u}_3 + \varepsilon_3 \omega_3;
\end{array}
\right.
\end{eqnarray*}
where
$\widetilde{u}_1=4.010906415,$ 
$\widetilde{u}_2=1.178843705,$
$\widetilde{u}_3=4.810249881$
(that is the stable  equilibrium state $E_2$),
$\varepsilon_1=\varepsilon_2=\varepsilon_3=0.001$,
and $\omega_i\in[0,1]$ are random variables, with $i=1,2,3$. 

First we consider the model with constant diffusions,
with $a_1=0.1$,   $a_2=0.0001$,  $a_3=0.1$ (see Figure \ref{fig:SARS_const_diff}).
We represent the random initial condition of the susceptible population
on the left top of the Figure  \ref{fig:SARS_const_diff}
(the graph of random initial condition is is indistinguishable to the three populations).
In the same figure, we can see the graph of the three population at time $t=2.5$.

On the other hand, in figure \ref{fig:SARS_inv_prop_diff} we observe another simulations with 
a nonlocal diffusion. 
In this simulation we consider the diffusion rate coefficients given by
\begin{equation}
\label{trunc}
\widetilde{a}_i(s) = 
\left\{
\begin{array}{ll}
M& \hbox{ if } s>M\\
s& \hbox{ if } \varepsilon\leqslant \dfrac{d_i}{(s-\widetilde{u}_i)^2} \leqslant M\\
\varepsilon& \hbox{ if } s<\varepsilon,
\end{array}
\right.
\end{equation}
with $M=10^4$, $\varepsilon=10^{-4}$,
$d_1=d_3=400000$ and $d_2=400$.
The choice of $d_i$ is in order to take diffusion coefficients close to the values 
of the constant coefficient diffusion of
the first simulation (Figure \ref{fig:SARS_const_diff}) at time $t=0$.
In this case of nonlocal diffusion, the coefficient diffusions are not constant in time,
obtaining obviously different results than the first simulation of this example 2.
Moreover, we remark very different behaviours between both simulations (compare
Figure  \ref{fig:SARS_const_diff} and Figure  \ref{fig:SARS_inv_prop_diff}).
Outside the nonlocal diffusion, we consider the same parameters for both simulations
in this example 2.
Finally, we observe in Figure \ref{fig:diff_ex2}, the
 evolution in time of the diffusion
$\displaystyle a_i\left(\int_\Omega u_i\,dx\right)$, with $i=1,2,3$ 
for the simulation with the nonlocal diffusion.

\renewcommand{\theequation}{A.\arabic{equation}}
\setcounter{equation}{0}
\section*{Appendix A. Uniqueness of weak solutions.} \label{uniq}

In this appendix we prove uniqueness of weak solutions
to our systems by using duality technique (see e.g. \cite{Rula}),
thereby completing the well-posedness
analysis. 

First, we consider $(u_{1},u_{2},u_{3})$ and $(v_{1},v_{2},v_{3})$ two solutions of the
system  (\ref {S1})-(\ref {S3}). We set $U_i=u_{i}-v_{i}$ for $i=1,2,3$, then $U_i$ satisfies
\begin{equation}\label{eq1:uniq-gener}
 \left\{\begin{array}{lll}\displaystyle
\displaystyle 
	&\displaystyle \pt U_1 
		-\Bigl(a_1\Bigl(\int_\Om u_1\dx\Bigl) \Delta u_1
		-a_1\Bigl(\int_\Om v_1\dx\Bigl) \Delta v_1\Bigl)\\ 
		&\qquad \qquad = 
		-(\sigma(u_1,u_2,u_3)-\sigma(v_1,v_2,v_3))-\mu U_1,\\ 
		&\displaystyle \pt U_2 
		-\Bigl(a_2\Bigl(\int_\Om u_2\dx\Bigl) \Delta u_2
		-a_2\Bigl(\int_\Om v_2\dx\Bigl) \Delta v_2\Bigl)\\ 
		&\qquad \qquad = 
		(\sigma(u_1,u_2,u_3)-\sigma(v_1,v_2,v_3))-(\gamma+\mu) U_2,\\ 
		&\displaystyle \pt U_3 
		-\Bigl(a_3\Bigl(\int_\Om u_3\dx\Bigl) \Delta u_3
		-a_3\Bigl(\int_\Om v_3\dx\Bigl) \Delta v_3\Bigl)=\gamma U_2,\\ 
		&\nabla u_i\cdot\eta=\nabla v_i \cdot
		\eta=0 \text{ on }\Sigma_T, \text{ $i=1,2,3$},\\
& U_i(x,0)=0 \mbox{ for } x \in \Omega,
		\text{ $i=1,2,3$}.
\end{array}\right.
\end{equation}
Now, we define the function $\varphi_i$ solution of the variational problem
\begin{equation}\label{ellip:uniq-gener}
\left\{
\begin{array}{l}\displaystyle 
\int_\Omega\nabla \varphi_i(t,\cdot) \cdot \nabla \phi\, dx =\int_\Omega U_i(t,\cdot)\phi \, dx,\\
\displaystyle 
\mbox{ for all $\phi\in H^1(\Omega)$, such that $\int_\Omega \phi\, dx =0$}
\end{array}
\right. 
i=1,2,3,
\end{equation}
for a.e. $t \in (0,T)$. Since $u_{i}$ and $v_{i}$ are in $L^2(Q_T)$, then we get from the theory
of linear elliptic equations, the existence, uniqueness and regularity of
solution $\varphi_i$ satisfying 
$$
\varphi_i \in C([0,T];H^1(\Omega)) \mbox{ with } \int_\Omega \varphi_i(t,\cdot) \,dx =0,
\mbox{ for $i=1,2,3$}.
$$
Note that from the boundary condition of $\varphi_i$ in (\ref
{ellip:uniq-gener}) and $U_i(0,\cdot)=0$ we deduce that
\begin{equation}\label{ellip1:uniq-gener}
\nabla \varphi_i(0,\cdot)=0 \mbox{ in }L^2(\Omega) \mbox{ for }i=1,2,3.
\end{equation}

Multiplying the first, second and third equations in (\ref {eq1:uniq-gener}) by $\psi_1,\psi_2,\psi_3
\in L^2(0,T;H^1(\Omega))$, respectively, and integrating over $Q_t:=(0,t)\times
\Omega$, we get
\begin{equation}\label{est1:uniq-gener}
  \begin{split}&\displaystyle
    \si \int_0^t \langle\partial_s  U_i,\psi_i \rangle \,ds 
    +\iint_{Q_t}\Bigl(a_1\Bigl(\int_\Om u_1\dx\Bigl) \Grad u_1
    -a_1\Bigl(\int_\Om v_1\dx\Bigl) \Grad v_1\Bigl)\cdot \nabla \psi_1  \,dx\,ds \\ 
    &\qquad +\iint_{Q_t}\Bigl(a_2\Bigl(\int_\Om u_2\dx\Bigl) \Grad u_2
    -a_2\Bigl(\int_\Om v_2\dx\Bigl) \Grad v_2\Bigl)\cdot \nabla \psi_2  \,dx\,ds \\ 
    &\qquad \qquad +\iint_{Q_t}\Bigl(a_3\Bigl(\int_\Om u_3\dx\Bigl) \Grad u_3
    -a_3\Bigl(\int_\Om v_3\dx\Bigl) \Grad v_3\Bigl)\cdot \nabla \psi_3  \,dx\,ds \\ 
    &\quad = 
    -\iint_{Q_t} \Bigl((\sigma(u_1,u_2,u_3)-\sigma(v_1,v_2,v_3))-\mu U_1\Bigl)\psi_1\,dx\,ds,\\ 
    &\qquad \quad + 
    \iint_{Q_t} \Bigl((\sigma(u_1,u_2,u_3)-\sigma(v_1,v_2,v_3))-(\gamma+\mu) U_2\Bigl)\psi_2\,dx\,ds 
    +\iint_{Q_t}\gamma U_2\psi_3\,dx\,ds.
  \end{split}
\end{equation}
Since $\varphi_i \in L^2(0,T;H^1(\Omega))$ we can take $\psi_i=\varphi_i$
in (\ref {est1:uniq-gener}) and we obtain from \eqref{ellip:uniq-gener} and
(\ref{ellip1:uniq-gener})
\begin{equation}\label{est2:uniq-gener}
\begin{array}{lll}\displaystyle
2\int_0^t \langle\partial_s  U_i,\varphi_i \rangle \,ds
&\displaystyle=-2\int_0^t \langle\partial_s  \Delta \varphi_i,\varphi_i \rangle \,ds \\
&\displaystyle=\int_\Omega |\nabla \varphi_j(t,x)|^2 \,dx 
-\int_\Omega \abs{\nabla \varphi_i(0,x)}^2 \,dx           \\
&=\displaystyle\int_\Omega |\nabla \varphi_i(t,x)|^2 \,dx.
  \end{array}
\end{equation}
From definition of $\sigma$ we obtain easily
\begin{equation}\label{est-sigma-hold}
\abs{\sigma(u_1,u_2,u_3)-\sigma(v_1,v_2,v_3)}\leq C\si \abs{u_i-v_i},
\end{equation}
for some constant $C>0$.

Using (\ref {ellip:uniq-gener}), (\ref {est-sigma-hold}),
H\"{o}lder's, Young's, Sobolev poincar\'e's inequalities yields from
(\ref{est1:uniq-gener}) with $\psi_i=\varphi_i$
\begin{equation}\label{est4:uniq-gener}
  \begin{split}\displaystyle
\displaystyle\si \int_0^t \langle\partial_s  U_i,\varphi_i \rangle \,ds
=&-\si \iint_{Q_t}a_i\Bigl(\int_\Om u_i\dx\Bigl) U_i\Delta  \varphi_i\,dx\,ds\\
& 
\displaystyle
-\si \iint_{Q_t}\Bigl(a_i\Bigl(\int_\Om u_i\dx\Bigl)-a_i\Bigl(\int_\Om v_i\dx\Bigl) \Grad v_i\cdot
\nabla \varphi_i\,dx\,ds
\\
&
-\iint_{Q_t} \Bigl((\sigma(u_1,u_2,u_3)-\sigma(v_1,v_2,v_3))-\mu U_1\Bigl)\varphi_1\,dx\,ds\\ 
    &+ 
    \iint_{Q_t} \Bigl((\sigma(u_1,u_2,u_3)-\sigma(v_1,v_2,v_3))-(\gamma+\mu) U_2\Bigl)\varphi_2\,dx\,ds 
    \\ 
    &+\iint_{Q_t}\gamma U_2\varphi_3\,dx\,ds\\ 
    =&-\si M_i\iint_{Q_t}\abs{U_i}^2\,dx\,ds\\
& 
\displaystyle
+\si \frac{M_i}{12}\iint_{Q_t}\abs{U_i}^2\,dx\,ds+C_1\si \int_{0}^t\norm{\Grad v_i}^2_{L^2(\Om)}
\norm{\Grad \varphi_i}^2_{L^2(\Om)}\,ds
\\
&+\Bigl(\si \frac{M_i}{12}\iint_{Q_t}\abs{U_i}^2\,dx\,ds
+\frac{M_1}{12}\iint_{Q_t}\abs{U_1}^2\,dx\,ds\Bigl)
\\&
+C_2\int_{0}^t
\norm{\Grad \varphi_1}^2_{L^2(\Om)}\,ds+\si \frac{M_i}{12}\iint_{Q_t}\abs{U_i}^2\,dx\,ds
+C_3\int_{0}^t
\norm{\Grad \varphi_2}^2_{L^2(\Om)}\,ds\\ 
&
+\frac{M_2}{12}\iint_{Q_t}\abs{U_2}^2\,dx\,ds
+C_4\int_{0}^t
\norm{\Grad \varphi_3}^2_{L^2(\Om)}\,ds\\ 
\leq& (C_1+C_2+C_3+C_4)\si \int_{0}^t(\norm{\Grad v_i}^2_{L^2(\Om)}+1)
\norm{\Grad \varphi_i}^2_{L^2(\Om)}\,ds,
  \end{split}
\end{equation}
for some constants $C_1,C_2,C_3,C_4>0$. 

Using $\nabla v_i \in L^2(Q_T)$ for $i=1,2,3$ and
Gronwall's lemma to conclude from (\ref{est4:uniq-gener})
$$
\varphi_i=0, \mbox{ $i=1,2,3$},
$$
almost everywhere in
$Q_T$, ensuring the uniqueness of weak solutions.

\subsection*{Acknowledgment}

This work  has been supported by  Fondecyt   projects \# 1070682
\# 1070694,
Fondap in Applied Mathematics (Project \# 15000001)
and INRIA/Conicyt Bendahmane-Perthame.

\def\cprime{$'$} \def\cprime{$'$} \def\cprime{$'$}

\end{document}